\documentclass[11pt]{article}

\usepackage{amsfonts}
\usepackage{graphicx, amssymb, latexsym,cite}

\begin{document}

\parskip 1mm 

 \def\KillFactor{{\sc KillFactor}}
\def\a{\alpha} 
 \def\b{\beta}
 \def\e{\epsilon}
 \def\d{\delta}
 \def\D{\Delta}
 \def\c{\chi}
 \def\k{\kappa}
 \def\g{\gamma}
 \def\t{\tau}
\def\ti{\tilde}
 \def\N{\mathbb N}
 \def\Q{\mathbb Q}
 \def\Z{\mathbb Z}
 \def\C{\mathbb C}
 \def\F{\mathbb F}
 \def\G{\Gamma}
 \def\go{\rightarrow}
 \def\do{\downarrow}
 \def\ra{\rangle}
 \def\la{\langle}
 \def\fix{{\rm fix}}
 \def\ind{{\rm ind}}
 \def\rfix{{\rm rfix}}
 \def\diam{{\rm diam}}
 \def\M{{\cal M}}
 
 \def\rank{{\rm rank}}
 \def\soc{{\rm soc}}
 \def\Cl{{\rm Cl}}
 \def\A{{\rm Alt}}
 \def\SL{{\rm SL}}
 \def\GL{{\rm GL}}
 
 \def\Sym{{\rm Sym}}
 \def\Ext{{\rm Ext}}
 \def\E{{\cal E}}
 \def\l{\lambda}
 \def\Lie{\rm Lie}
 \def\s{\sigma}
 \def\O{\Omega}
 \def\o{\omega}
 \def\ot{\otimes}
 \def\op{\oplus}
 \def\pf{\noindent {\bf Proof.$\;$ }}
 \def\Proof{{\it Proof. }$\;\;$}
 \def\no{\noindent}
\def\hal{\unskip\nobreak\hfil\penalty50\hskip10pt\hbox{}\nobreak
 \hfill\vrule height 5pt width 6pt depth 1pt\par\vskip 2mm}

 \renewcommand{\thefootnote}{}

\newcommand{\Oh}{O}

\newtheorem{theorem}{Theorem}
\newtheorem{remark}[theorem]{Remark}
 \newtheorem{thm}{Theorem}[section]
 \newtheorem{theor}{Theorem}[subsection]
 \newtheorem{prop}[thm]{Proposition}
 \newtheorem{lem}[theor]{Lemma}
 \newtheorem{lemma}[thm]{Lemma}
 \newtheorem{propn}[theor]{Proposition}
 \newtheorem{cor}[thm]{Corollary}
 \newtheorem{coroll}[theorem]{Corollary}
 \newtheorem{coro}[theor]{Corollary}
 \newtheorem{rem}[thm]{Remark}
 \newtheorem{cla}[thm]{Claim}
 
\title{Recognition of finite exceptional groups of Lie type}
\author{Martin W. Liebeck \\ Department of Mathematics\\
 Imperial College \\ London SW7 2BZ \\ UK
 \and
 E.A. O'Brien \\
 University of Auckland \\
 Auckland \\ New Zealand}
\date{}
\maketitle

\begin{abstract}
Let $q$ be a prime power and let $G$ be an absolutely irreducible subgroup 
of $GL_d(F)$, where $F$
is a finite field of the same characteristic as $\F_q$, 
the field of $q$ elements. 
Assume that $G \cong G(q)$, a quasisimple group of exceptional Lie type
over $\F_q$ which is neither a Suzuki nor a Ree group. 
We present a Las Vegas algorithm that constructs an
isomorphism from $G$ to the standard copy of $G(q)$.
If 
$G \not\cong {}^3\!D_4(q)$ with $q$ even, 
then the algorithm runs in polynomial time,
subject to the existence of a discrete log oracle.
\end{abstract}


\section{Introduction}\label{intro}
Informally, a {\it constructive recognition algorithm}
constructs an explicit isomorphism
between a quasisimple group $G$ and a `standard' 
copy of $G$, and exploits this isomorphism
to write an arbitrary element of $G$ as
a word in its defining generators.
For a more formal definition, see \cite[p.\ 192]{Seress03}.
Such algorithms play a critical role in 
the `matrix group recognition project' which aims to
develop efficient algorithms for the investigation of 
subgroups of $\GL_d(F)$ where $F$ is a finite field. 
We refer to the recent survey \cite{OBrien11}
for background related to this work.
Such algorithms are available for classical groups; see, for example,
\cite{DLGO, DLGO14, LGO}. 
Here we present constructive recognition algorithms 
for the finite exceptional groups of Lie type. 

Let $G(q)$ denote a quasisimple exceptional group of Lie type over $\F_q$,
a finite field of size $q$. 
Howlett, Rylands \& Taylor  \cite{HRT} provide defining matrices
for a specific faithful irreducible representation of minimal dimension 
of the simply connected group of type $G(q)$: we call this
representation the {\it standard copy} of type $G(q)$.

Our principal result is the following. In the statement, $V_d(F)$ denotes the 
underlying vector space of dimension $d$ over the field $F$.

\begin{theorem}\label{main}
Let $q$ be a prime power and let $G$ be an absolutely irreducible subgroup of 
$GL_d(F)$, where $F$
is a finite field of the same characteristic as $\F_q$. Assume
that $G \cong G(q)$, a quasisimple group of exceptional Lie type
over $\F_q$ for $q > 2$, excluding Suzuki and Ree groups, 
and also $^3\!D_4(q)$ with $q$ even. 
There is a Las Vegas algorithm that constructs an
isomorphism from $G$ to the standard copy of type $G(q)$ modulo a 
central subgroup, and also constructs the inverse isomorphism; it also
computes the high weight of $V_d(F)$ as a $G$-module, 
up to a possible twist by a field or graph automorphism. 
The algorithm runs in polynomial time,
subject to the existence of a discrete log oracle for extensions of 
$\F_q$ of degree at most $3$. 
\end{theorem}

The possible central subgroups of the standard copy of type $G(q)$ are 
trivial except when $G(q)$ is of type
$E_6^\e(q)$ ($\e = \pm 1$) or $E_7(q)$, in which case they have order 
dividing $(3,q-\e)$ and $(2,q-1)$ respectively.

We now discuss how the isomorphisms in the statement of the theorem 
are realised. 
Let $\hat G(q)$ denote the standard copy of type $G(q)$. 
It has a {\it Curtis-Steinberg-Tits} presentation,  
which involves only those relations which arise from certain 
rank 2 subgroups of $G(q)$: 
namely, the commutator relations among root elements corresponding to pairs of 
fundamental roots in the corresponding Dynkin diagram. 
Babai {\it et al.}\cite[\S 4.2 and 6.1]{BGKLP} reduce this
presentation by running over root elements parametrised by an 
$\F_p$-basis for $\F_q$ (where $p$ is the characteristic of $\F_q$). 
Those root elements of $\hat G(q)$ which satisfy this 
reduced Curtis-Steinberg-Tits presentation are the {\it standard generators} 
$\hat {\cal S}$ for $\hat G(q)$. 

Given a group $G$ as in the statement of the theorem, 
described by a generating set $X$, 
our algorithm produces a collection ${\cal S}$ of generators of $G$ 
(as words in $X$) which 
satisfy the reduced Curtis-Steinberg-Tits presentation. 
These are then used to construct 
the required isomorphisms $\phi: G \rightarrow \hat G(q)/Z$ and 
$\psi: \hat G(q)/Z \rightarrow G$ (where $Z$ is a central subgroup), as follows. 
Cohen, Murray \& Taylor \cite{CMT} developed the
{\it generalised row and column reduction algorithm}: 
in polynomial time, for a given high weight representation of 
$G \cong G(q)$ with $G(q)$ of untwisted type, 
this algorithm writes an arbitrary $g \in G$ as a word $w({\cal S})$ in  
the standard generators; this has now 
been extended to twisted types in \cite{CT}.
Now $\phi(g) = w(\hat {\cal S})Z$, the corresponding word in the 
standard generators of 
$\hat G(q)$, defines the isomorphism $\phi$. 
The inverse isomorphism $\psi$ is defined 
similarly: for $\hat g \in \hat G(q)$, the algorithms of \cite{CMT, CT} 
express $\hat g$ as a 
word $w(\hat {\cal S})$, and we set $\psi(\hat g Z) = w({\cal S})$.

Together, the algorithms of Theorem \ref{main} and of \cite{CMT, CT} provide 
a solution 
to the {\it constructive membership problem} for $G = \langle X \rangle$:   
namely, express an arbitrary $g\in G$ as a word in~$X$.

Our algorithms to find standard generators in $G$ begin by constructing 
$SL_2$ subgroups of $G$ which can be placed as nodes in
the Dynkin diagram so that they pairwise generate the appropriate
group of rank 2, and these are then used to label root elements and
toral elements of $G$ relative to a fixed root system. 
We use the root elements to compute the
high weight of the given representation of $G$ on $V_d(F)$, 
and then exploit the algorithms of \cite{CMT, CT} to 
set up the isomorphisms  explicitly.

To construct the $SL_2$ subgroups and label root
elements, we use involution centralizers in $G$. 
That such centralizers can be constructed in 
Monte Carlo polynomial time follows in odd characteristic from \cite{PW}, 
and in even characteristic from \cite{LO}; 
see Section \ref{centralizer} for further discussion.

A distinguishing feature of our work is that the resulting 
algorithms are practical; this desire significantly influenced our design. 
Our algorithms are implemented and will be publicly available 
in {\sc Magma} \cite{Magma}. 

The excluded Suzuki and Ree groups of types 
$^2\!B_2(q)$, $^2\!G_2(q)$, and $^2\!F_4(q)$ 
were studied by B\"a\"arnhielm \cite{HBthesis, HB06, HB12}. 
His constructive recognition algorithms apply to conjugates 
of the standard copy of $^2\!B_2(q)$ and $^2\!G_2(q)$,
and run in polynomial time subject to the availability 
of a discrete log oracle. For the groups $^3\!D_4(q)$ ($q$ even), also excluded
in the theorem, we provide a practical algorithm with running time $O(q)$.
We also present practical algorithms for groups defined
over $\F_2$, the only field not covered by Theorem \ref{main}. 
Where feasible, our theoretical results also include $\F_2$. 

As stated, Theorem \ref{main} applies to absolutely irreducible
representations of quasisimple groups of exceptional Lie type.
The principal motivation for stating it under this assumption is our 
application of the algorithms of \cite{CMT, CT} to realise the isomorphisms 
between $G$ and $\hat{G}(q)/Z$.
Using the Meataxe and associated machinery \cite[Chapter 7]{HoltEickOBrien05}, 
the result can easily 
be reformulated to apply, with unchanged complexity, to all matrix 
representations
(not necessarily irreducible)  
in defining characteristic. For all but $E_8(q)$ in even 
characteristic, 
our algorithms to construct the $SL_2$ subgroups and to label 
the root and toral elements are black-box {\it provided} that the 
algorithms employed
in Theorem \ref{available-crec} for constructive recognition of 
small rank classical groups are black-box. 
Since algorithms are available for these tasks
(see, for example, \cite{DLGO14} and its references), 
a version of Theorem \ref{main} could be formulated for black-box groups. 
We refrain from doing so.

Kantor \& Magaard \cite{KM} presented 
black-box Las Vegas algorithms to recognise constructively
the exceptional simple groups of Lie type and rank at least 2,
other than ${}^2F_4(q)$, defined over a field of known size.
These have complexity depending linearly on the size of the field. 
Dick \cite{dick} developed a polynomial-time algorithm, 
a modification of that proposed by \cite{KM}, 
for $F_4(q)$ in odd characteristic. 
Relying as it does on centralizers of involutions, 
our work differs substantially from \cite{KM}.

The structure of the paper is as follows.
Section \ref{prelim} records a number of 
results which underpin our algorithm. 
In Section \ref{generate} we prove 
results on probabilistic generation for certain
groups of Lie type. In Sections \ref{sl3sl6}-\ref{3d4sec} we present
algorithms to construct $SL_2$ subgroups of a group $G$ 
as in Theorem \ref{main} 
which correspond to the nodes in the associated Dynkin diagram. 
Sections \ref{lab}-\ref{pres} contain algorithms to 
label root elements and toral elements of $G$ relative to 
a fixed root system; to determine the
high weight of the given representation of $G$ on $V_d(F)$; 
and to construct the standard generators for $G$.
In Section \ref{qeq2} we present algorithms for the special case 
of groups defined over $\F_2$.
In Section \ref{imp} we report on our 
implementation in {\sc Magma}. Finally, 
for each group of exceptional Lie type and Lie rank at least 2, 
its reduced Curtis-Steinberg-Tits presentation
on standard generators is listed explicitly in Appendix~\ref{app}.

\vspace{4mm}
\no {\bf Acknowledgements.} 
This work was supported in part by the Marsden Fund of
New Zealand via grant UOA 1323 and an International Exchange grant 
from the Royal Society. We thank Arjeh Cohen and Don Taylor 
for their work on \cite{CT}, providing access to
their implementation of the algorithms of \cite{CMT, CT}, and 
for many helpful discussions; we thank 
Scott Murray for assistance on Lie algebra computations.

\section{Background and preliminaries} \label{prelim}

A {\it Monte Carlo} algorithm is a randomised algorithm which always
terminates but may return a wrong answer with probability less than
any specified value.  A {\it Las Vegas} algorithm is a randomised algorithm
which never returns an incorrect answer, but may report
failure with probability less than any specified value.

Our algorithms usually search for elements of $G$ having
a specified property. If $1/k$ is a lower  bound
for the proportion of such elements in $G$, then we can readily
prescribe the probability of failure of the corresponding algorithm.
Namely, to find such an element by random
search with a probability of failure less
than a given $\epsilon \in (0, 1)$
it suffices to choose (with replacement) a sample
of uniformly distributed random elements in $G$
of size at least $\lceil -\log_e(\epsilon)k \rceil$.

Babai \& Szemer{\'e}di \cite{BabaiSzemeredi84} introduced
the {\it black-box group} model, where
group elements are represented by bit-strings
of uniform length.
The only group operations
permissible are multiplication, inversion, and checking
for equality with the identity element.
Seress \cite[p.\ 17]{Seress03} defined
a {\it black-box algorithm} as one which does not
use specific features of the group representation,
nor particulars of how group operations
are performed; it can only use the operations
listed above.  

Babai \cite{Babai91} present a Monte Carlo algorithm to
construct in polynomial time independent nearly uniformly distributed
random elements of a finite group.  An alternative is the
{\it product replacement algorithm} of Celler
{\it et al.\ }\cite{Celleretal95}, which runs in  polynomial time
by a result of \cite{Pak00}. For a discussion of both algorithms
we refer to \cite[pp.\ 26--30]{Seress03}.

Often it is necessary to investigate the order of $g\in \GL_d(\F_q)$,
which, due to integer factorisation,  cannot be
determined in polynomial time. We can, however, determine its
\emph{pseudo-order}, a good multiplicative upper bound for $|g|$,
and the exact power of any specified prime
that divides $|g|$, using a Las Vegas algorithm with complexity
$\Oh(d^3 \log d + d^2 \log d \log q)$. 
Our results sometimes assume the existence of an {\it order oracle} 
but, in our applications, it always suffices to use pseudo-order. 
A Las Vegas algorithm with the
same complexity allows us to compute large powers $g^n$ where $0\leq n<q^d$.
We refer to  \cite[\S 2 and 10]{LGO} for more details and references.

Leedham-Green \& O'Brien \cite{LeedhamGreenOBrien02} present
Monte Carlo algorithms to generate random elements of the normal
closure of a subgroup, and to determine membership in a normal
subgroup of a black-box group having an order oracle.
Babai \& Shalev \cite{BabaiShalev01} prove that if the 
normal subgroup is simple
and non-abelian, then the membership algorithm
runs in Monte Carlo polynomial time.
A consequence is a Monte Carlo black-box
algorithm to prove that a group is perfect.
This algorithm is used together with the black-box polynomial-time 
algorithm described in \cite[pp.\ 38-40]{Seress03}
to construct the derived series of a group. 

To construct a direct factor of a semisimple group,
we use the black-box algorithm, \KillFactor,
of \cite[Claim 5.3]{BabaiBeals99}; 
that it runs in polynomial time is a consequence
of \cite[Corollary 4.2]{Babai09}.

If a matrix group acts absolutely irreducibly  on its underlying
vector space of dimension $d$, then we can determine the classical forms
it preserves in $\Oh(d^3)$ field operations
(see \cite[\S 7.5.4]{HoltEickOBrien05}).
A hyperbolic basis for a vector space of dimension $d$ with a given
non-degenerate bilinear form can be constructed
in $\Oh(d^3)$ field operations (see \cite{Brooksbank03}
for an algorithm to perform this task).

\subsection{Recognition for classical groups}\label{classical-crec}
Babai {\it et al.\ }\cite{BabaiKantoretal02} proved the following.
\begin{thm}\label{identify}
Given a black-box group $G$ isomorphic
to a simple group of Lie type of known characteristic,
the standard name of $G$ can be computed using a 
Monte Carlo polynomial-time algorithm.
\end{thm}

Liebeck \& O'Brien \cite{LO} present a Monte Carlo
black-box polynomial-time algorithm to identify the defining characteristic.
Kantor \& Seress \cite{KS09} give  
an alternative algorithm for absolutely irreducible matrix groups.

Kantor \& Seress \cite{KS01} developed 
the first black-box Las Vegas algorithms to recognise constructively
classical groups;  these have complexity depending linearly on the 
size of the field. 
More recently, Leedham-Green \& O'Brien \cite{LGO} developed 
algorithms for classical groups in natural representation and odd 
characteristic; those of Dietrich {\it et al.\ }\cite{DLGO}  
apply to even characteristic.
Black-box equivalents appear in \cite{DLGO14}.
All run in time polynomial in the size of the input subject
to the availability of a discrete log oracle. 

Our algorithms for the labelling of root and toral elements 
rely on the availability of constructive recognition algorithms 
for the classical groups listed in the following theorem.
As defined in \cite{LGO}, the {\it standard copy} of a 
classical group is its natural matrix representation,
preserving a specified form.

\begin{thm} \label{available-crec}
Let $q$ be a prime power, and let $G$ be a subgroup 
of $GL_d(F)$, where $F$ is a field of the same characteristic as $\F_q$. 
Assume that $G$ is isomorphic to one of the following classical groups:
$SL_2(q)$, $SL_3(q)$, $Sp_4(q)$, or $SU_4(q)$ in all characteristics;
$SL_6(q)$, $Sp_6(q)$, $SU_6(q)$, 
$\O_{8}^{\pm}(q)$ or $\O_{12}^+(q)$ for even $q$.
There is a Las Vegas algorithm that constructs an
isomorphism from $G$ to its standard copy. 
Subject to the existence of a discrete log oracle, 
the algorithm runs in polynomial time.
\end{thm}

This follows from 
\cite{Brooksbank03a}, 
\cite{BrooksbankKantor06},
\cite{Brooksbank08}, 
\cite{CLGO}, 
\cite{DLGO14}, 
and \cite{sl3q}. 

\subsection{Groups of Lie type}\label{notnlie}

We use $SL_n^\e(q)$ to denote $SL_n(q)$ for 
$\e=1$ and $SU_n(q)$ for $\e=-1$; 
we adopt similar conventions for $D_4(q)$ and $^2\!D_4(q)$;  
and  for $E_6(q)$ and $^2\!E_6(q)$.
Dynkin diagrams for exceptional root systems are labelled as follows:
\[
\begin{array}{llllllllcllllll}
E_l &&&&1& - &3& - &4& - & 5& - & \cdots & - &l \\
&&&& && && | & & & & \\
&&&& &&& & 2 & &&& 
\end{array}
\]
\[
\begin{array}{llllllllllllllll}
F_4&&&&1 & - &2 & =>= & 3& - & 4&&&&&
\end{array}
\]
\vspace{4mm}
\[
\begin{array}{lllllllllllllllllll}
G_2&&&&& 1 & \equiv > \equiv & 2 &&&&&&&&&&
\end{array}
\]
where each node $i$ represents a simple root $\a_i$. This is the labelling of
Bourbaki \cite[p.~250]{bour}, except for 
$G_2$, where $\a_1$ and $\a_2$ are interchanged.

Let $G = G(q)$ be an exceptional group of Lie type over $\F_q$; 
we exclude  Suzuki and Ree groups. The root system of $G(q)$ is 
described in \cite[Chapter 3]{car};
if $G(q)$ is of twisted type, then we use the twisted root system 
of \cite[Chapter 13]{car}.
For a long root $\a$ in the root system, we denote by 
$U_\a$ the corresponding long root group, and 
a conjugate of $\la U_{\pm \a} \ra \cong SL_2(q)$ is 
a {\it long} $SL_2$ subgroup of $G(q)$. 
For a fixed isomorphism between $\la U_{\pm \a} \ra$ and $SL_2(q)$,
we denote by $h_\a(c)$ the element of $\la U_{\pm \a}\ra$ corresponding 
to the matrix ${\rm diag}(c^{-1},c)$; 
if $\a$ is a fundamental root $\a_i$, then we may write 
$h_i(c)$. If $\a_1\ldots ,\a_l$ are 
fundamental roots and $c_1,\ldots ,c_l$ are integers, 
then $h_{c_1\ldots c_l}(\l) := h_1(\l^{c_1})\cdots h_l(\l^{c_l})$.

An involution in a long $SL_2$ subgroup of $G(q)$ is a 
{\it root} involution. These involutions and their 
centralizers play a major role in our work. 
Proposition \ref{invcent} lists the root involution centralizers; 
it appears in \cite[4.5]{GLS} (for $q$ odd) and \cite{AS} (for $q$ even). 
A {\it subsystem} subgroup is one generated by root groups 
corresponding to roots in a closed subsystem of the root system of $G(q)$.

\begin{prop}\label{invcent}
Let $G = G(q)$ be an exceptional group of Lie type over $\F_q$, 
and let $t$ be a root involution. Let $D$ be a subsystem 
subgroup of $G$ as in the following table:
\[
\begin{array}{|l|l|l|l|l|l|l|}
\hline
G & G_2(q) & ^3\!D_4(q) & F_4(q) & E_6^\e(q) & E_7(q) & E_8(q) \\
\hline
D & A_1(q) & A_1(q^3) & C_3(q) & A_5^\e(q) & D_6(q) & E_7(q) \\
\hline
\end{array}
\]

{\rm (i)} If $q$ is odd, then $C_G(t)$ has a subgroup $SL_2(q)D$ 
of index at most $2$; the factors $SL_2(q)$ and $D$ commute (elementwise).

{\rm (ii)} If $q$ is even, then $C_G(t) = QD$, 
where $Q$ is a normal $2$-subgroup of $C_G(t)$.
\end{prop}

If the root system of $G(q)$ has roots of different lengths, then 
a {\it short} $SL_2$ subgroup is one generated by a pair of 
opposite short root subgroups of $G(q)$; if $G(q)$ is of untwisted type, 
then these are isomorphic to $SL_2(q)$, otherwise they are 
isomorphic to $SL_2(q^2)$, 
or $SL_2(q^3)$ for $^3\!D_4(q)$.

Let $l$ be the rank of the root system of $G(q)$ and let $1,\ldots,l$ be the 
nodes of the Dynkin diagram.  Let 
$K_1,\ldots, K_l $ be long (short) $SL_2$ subgroups 
of $G(q)$ which satisfy the following:
\begin{enumerate}
\item $K_i$ is long (short) if and only if node $i$ is a long (short) root;

\item if nodes $i,j$ are not joined then $K_i$ and $K_j$ commute;

\item if nodes $i,j$ are joined then $\langle K_i,K_j\rangle$ is the 
appropriate rank 2 group of Lie type: $A_2(q)$ or $A_2(q^2)$ if 
$i,j$ are joined by a single bond; $B_2(q)$ or $^2\!A_3(q)$ if joined 
by a double bond; $G_2(q)$ or $^3\!D_4(q)$ if joined by a triple bond.
\end{enumerate}
We call such $K_1,\ldots, K_l$ {\it basic} $SL_2$ subgroups of $G(q)$. 

\subsection{Centralizers of involutions}\label{centralizer}
The centralizer of an involution in a black-box group having 
an order oracle can be constructed using an algorithm of 
Bray \cite{Bray}; he proved the following.

\begin{thm}
\label{thm:bray}
If $x$ is an involution in a group $H$, and $w$ is an arbitrary element of $H$,
then $[x,w]$ either has odd order $2k+1$, in which case
$w[x,w]^k$ commutes with $x$, or has even order $2k$, in which case
both $[x,w]^k$ and $[x,w^{-1}]^k$ commute with $x$.
If $w$ is uniformly
distributed among the elements of the group for which $[x,w]$
has odd order, then $w[x,w]^k$ is uniformly distributed among the
elements of the centralizer of $x$.
\end{thm}

Thus if the odd order case occurs sufficiently often, 
then we can construct random elements of 
the involution centralizer in Monte Carlo polynomial time.

Parker \& Wilson \cite[Theorems 1-4]{PW} proved the following two results.

\begin{thm}\label{clasthm}
There is an absolute constant $c>0$ such that if $H$ is a finite
simple group of Lie type of Lie rank $r$
defined over a field of odd characteristic,
and $x$ is an involution in $H$, 
then the proportion of $h \in H$ such that $[x,h]$ has 
odd order is at least $c/r$.
\end{thm}

\begin{thm}\label{powerup}
There is an absolute constant $c>0$ such that if $H$ is a finite
simple group of Lie type of Lie rank $r$
defined over a field of odd characteristic, and $C$ is a conjugacy class of 
involutions in $H$, 
then the proportion of elements of $H$ which power up to an element of 
$C$ is at least $c/r^{3}$.
\end{thm}

By Theorem \ref{powerup}, an involution in a specified class of $H$ 
can be constructed in polynomial 
time by searching for an element of even order and computing 
a suitable power. By Theorem \ref{clasthm}, random elements of the 
centralizer of this involution
can be constructed, and a bounded number of these generate the 
centralizer (see \cite{LSh}).

We shall make frequent use of the following lemma, also proved by Parker \& Wilson \cite[Lemma 26]{PW}.

\begin{lemma} \label{parwil}
Let $H$ be a finite group of Lie type, $T$ a maximal torus in $H$ and let $C$ be a conjugacy class of $H$. 
Assume that at
least a proportion $k$ of the regular semisimple elements of $T$ power to a member of $C$. Then at least a proportion
$k/|N_H(T): T|$ of the elements of $H$ power to an element of $C$.
\end{lemma}

Liebeck \& O'Brien \cite{LO} proved the following.

\begin{thm}\label{even-cent} If $H$ is a finite group of Lie type over a 
field of characteristic $2$, and $t$ is a root involution 
of $H$ such that $C_H(t)$ is not soluble,
then the proportion of $h \in H$ such that $[t,h]$ has 
odd order is at least $1/4$.
\end{thm}

Hence random elements of the centralizer of a root involution can 
be constructed in polynomial time, and a bounded
number of these will generate the centralizer (see \cite[Lemma 3.10]{LO}).
\subsection{The Formula}
Variations of the following lemma, sometimes known as 
the ``Formula", have been used in algorithms for some years -- 
see, for example, \cite[Section 4.10]{formula}. 
A proof can be found in \cite[7.1]{DLGO}.

\begin{lemma}\label{Formula}
Let $K = H\ltimes M$ where $M$ has exponent $2$. Suppose $h\in H$
has odd order and acts fixed point freely on $M$. If $k=am \in K$
where $a\in C_H(h)$ and $m\in M$, then $a=hk(hh^k)^{(|h|-1)/2}$.
\end{lemma}

The lemma sometimes allows us to construct a complement to a 
normal 2-subgroup in a semidirect product. 
We apply it when $H = \langle h\rangle \times D$ 
for quasisimple $D$. Now the lemma enables us to construct 
random elements of $H$ (namely, $hk(hh^k)^{(|h|-1)/2}$ for random $k\in K$). 
Since, by \cite{LSh}, $D$ may be generated by two random elements, we 
can thus construct a generating set for $D$. 

\section{Probabilistic generation of certain groups}\label{generate}
Our algorithms rely on various results on probabilistic
generation; these we now present. 

\begin{prop}\label{d4gen} 
Let $G=D_4^\e(q)$ with $\e=\pm$ and $q>2$ even, 
and let $x$ be an element of order $q+1$ in a long $SL_2$ 
subgroup of $G$. 
For random $g\in G$, the probability that 
$\la x,x^g\ra = G$ is positive, and is at least $1-c/q$, where 
$c$ is an absolute constant. 
\end{prop}

\pf 
Let $P$ be the probability
that $\la x,x^g\ra = G$ for random $g \in G$. If $\la
x,x^g\ra \ne G$, then $x, x^g \in M$ for some maximal subgroup $M$ of
$G$. Given a maximal subgroup $M$ containing $x$, the probability that
$x^g$ lies in $M$ is $|x^G\cap M|/|x^G|$. 
It is well-known (see \cite[2.5]{LSax}) that
\begin{equation}\label{wellknown}
\frac{|x^G\cap M|}{|x^G|} = \frac{{\rm fix}_{G/M}(x)}{|G:M|},
\end{equation}
 where ${\rm fix}_{G/M}(x)$ denotes the number of fixed points of $x$ in the
action of $G$ on the cosets of $M$. Also, the number of conjugates of
$M$ containing $x$ is ${\rm fix}_{G/M}(x)$. Hence, if ${\cal
 M}$ is a set of representatives of the conjugacy classes of maximal
subgroups of $G$, then
\begin{equation}\label{probbd}
1-P \le \sum_{M\in {\cal M}} \frac{{\rm fix}_{G/M}(x)^2}{|G:M|}.
\end{equation}
The maximal subgroups of $G$ are determined in \cite{Kl} (for $\e=+$)
and are listed in \cite[Tables 8.52--8.53]{BHRD} (for $\e = -$). 
In Tables \ref{fixplus} and \ref{fixminus}, we list those maximal 
subgroups $M$ which contain a conjugate of $x$, together with the values 
of ${\rm fix}_{G/M}(x)$ and 
$|G:M|$. The notation is standard: $P_i$ denotes a parabolic subgroup, 
the stabilizer of a totally singular $i$-space; 
and $N_i^\d$ ($\d=\pm$) is the stabilizer of a nonsingular subspace 
of dimension $i$ and type $\d$. Note that we omit $N_2^+$ from the tables: if 
$x,x^g$ are contained in $N_2^+$ then they lie in a subgroup of $N_2^+$ which 
is contained in $P_1$. In Table \ref{fixplus} we comment if 
a row covers 3 classes of maximal subgroups; in each case these 
are permuted by a triality automorphism of $G$.

\begin{table}[htb]
\[
\begin{array}{llll}
\hline 
M & {\rm fix}_{G/M}(x) & |G:M| & \hbox{Comment} \\
\hline
P_1 & (q+1)^2 & \frac{(q^4-1)(q^3+1)}{q-1} & 3 \hbox{ classes of  subgroups }M  \\
P_2 & 3(q+1) & \frac{(q^6-1)(q^2+1)^2}{q-1} & \\
N_1 & q(q^2-1) & q^3(q^4-1) &  3 \hbox{ classes of  subgroups }M   \\
N_2^- & \frac{1}{2}q(q-1)(q^2-q+2) & \frac{q^6(q^4-1)(q^3-1)}{2(q+1)} &  3 \hbox{ classes of  subgroups }M   \\
N_4^+.2 & \frac{1}{4}q^3(q-1)^3 & \frac{q^8(q^6-1)(q^2+1)^2}{4(q^2-1)} & \\
N_4^-.2 & \frac{1}{4}q^3(q+1)(q^2-1) & \frac{1}{4}q^8(q^6-1)(q^2-1) & 3 \hbox{ classes of  subgroups }M   \\
\hline
\end{array}
\]
\caption{Fixed points of $x$ for $G = \O_8^+(q)$} \label{fixplus}
\end{table}

\begin{table}[htb]
\[
\begin{array}{lll}
\hline 
M & {\rm fix}_{G/M}(x) & |G:M|  \\
\hline
P_1 & q^2+1 & \frac{(q^4+1)(q^3-1)}{q-1}  \\
P_2 & q+1 & \frac{(q^6-1)(q^4+1)}{q-1}  \\
P_3 & (q^2+1)(q+1) & (q^4+1)(q^3+1)(q^2+1) \\
N_1 & q(q^2+1) & q^3(q^4+1) \\
N_2^- & \frac{1}{2}q^2(q^2+1) & \frac{q^6(q^4+1)(q^3+1)}{2(q+1)}  \\
N_4^+ & \frac{1}{2}q^3(q^2+1)(q-1)+1 & \frac{q^8(q^6-1)(q^4+1)}{2(q^2-1)}  \\
\O_4^-(q^2).2 &  \frac{1}{2}q^3(q^2-1)(q+1) & \frac{1}{2}q^8(q^6-1)(q^2-1) \\
U_3(q) & 2(3,q+1)q^2(q^4-1)  &  (3,q+1)q^9(q^8-1)(q^3-1)  \\
\hline
\end{array}
\]
\caption{Fixed points of $x$ for $G = \O_8^-(q)$} \label{fixminus}
\end{table}

The values of ${\rm fix}_{G/M}(x)$ given in the tables are calculated 
reasonably routinely; we give a sketch. 
Let $G$ have natural module $V_8$. 
Regard $x$ as acting on an orthogonal decomposition 
$V_8 = V_4 + V_4'$, 
where $V_4,V_4'$ are non-degenerate subspaces of type $O_4^\e$ and 
$O_4^+$ respectively; 
$x$ acts trivially on $V_4$, and as an element in one of the $SL_2(q)$ factors 
of $\O_4^+(q) = SL_2(q) \otimes SL_2(q)$ on $V_4'$. 
For $M=P_1$ (or $N_1$), ${\rm fix}_{G/M}(x)$ is the number of singular 
(or nonsingular) 1-spaces in $V_4$; for $M=P_2$, 
${\rm fix}_{G/M}(x)$ is the sum of the numbers of singular 2-spaces fixed 
by $x$ in $V_4$ and $V_4'$; and for $M=P_3$, 
the singular 3-spaces in $V_8$ fixed by $x$ are spanned by one of the 
$q+1$ fixed 2-spaces in $V_4'$ together with a fixed 
1-space in $V_4$. For $M = N_2^-$, the 
2-spaces of type $O_2^-$ fixed by $x$ either lie in 
$V_4$ or in $V_4'$, and there are $q(q-1)$ of the latter, as can be seen 
using (\ref{wellknown}). Likewise, $N_4^\pm$-spaces in $V_8$ fixed by $x$ are 
sums of fixed nonsingular 2-spaces in $V_4$ and $V_4'$, and it is 
straight-forward to 
count these. Finally, the cases where 
$\e=-$ and $M = \O_4^-(q^2).2$ or $U_3(q)$ 
are handled using (\ref{wellknown}). 
In the first case, $x^G\cap M$ is a class of elements of order $q+1$ in 
$\O_4^-(q^2) \cong L_2(q^4)$, so has size 
$q^4(q^4+1)$. The second case arises from the adjoint representation 
of $U_3(q)$, in which $x$
acts as ${\rm diag}(\a,\a,\a^{-2})$ for some scalar $\a$ of order $q+1$ 
or its inverse. Hence 
$|x^G\cap M| = 2|SU_3(q):GU_2(q)| = 2q^2(q^2-q+1)$, 
from which ${\rm fix}_{G/M}(x)$ 
follows using (\ref{wellknown}).

The lower bound $1-c/q$ in the statement of the proposition follows from the 
information in the tables together with (\ref{probbd}). These also imply that 
$1-P$ is less than 1 for $q\ge 8$, giving the positivity statement for these 
values of $q$. For $q=4$ we can verify computationally that $G$ is 
generated by two conjugates of $x$. \hal

\vspace{4mm} For $q=2$ the probability in the previous proposition 
remains positive for $\e=-$, but is zero for $\e=+$.  

\begin{prop}\label{probp2} 
Let $G$ be one of $F_4(q), E_6^\e(q), E_7(q)$ or $E_8(q)$ with $q$ even, and
let $x$ be an element of order $q+1$ in a long $SL_2$
subgroup of $G$. For random $g\in G$, the probability that
$\la x,x^g\ra$ is a subsystem subgroup $D_4^\e (q)$ (for some $\e \in \{+,-\}$) is positive, and is at
least $1/6-c/q$, where $c$ is an absolute constant.
\end{prop}

\pf That the probability is positive follows immediately from 
Proposition \ref{d4gen}
(and the ensuing remark for the case $q=2$).
Let $D \cong D_4^\e (q)$ be a fixed subsystem
subgroup of $G$ which contains a long $SL_2$ subgroup containing $x$,
and define
\[
\D = \{ D^g : g \in C_G(x)\}.
\]
Observe that $|\D| = |C_G(x):C_G(x)\cap N_G(D)|$ 
since $C_G(x)$ acts transitively on $\D$.

We consider first the case where $\e = +$. Here $N_G(D) = D.S_3$,
$DT_2.S_3$, $DA_1(q)^3.S_3$ or $DD_4(q).S_3$ according as $G$ is of
type $F_4, E_6^\e, E_7$ or $E_8$, where $T_2$ denotes a rank 2 torus 
(see \cite[Table 5.1]{LSS}). 
The table below gives $C_G(x)$, $C_G(x) \cap N_G(D)$ and $|\D|$. 
\[
\begin{array}{|l|l|l|l|}
\hline 
G & C_G(x) & C_G(x) \cap N_G(D) & |\D|\sim \\
\hline 
F_4(q) & \la x \ra C_3(q) & (\la x \ra A_1(q)^3).S_3 & q^{12}/6 \\
E_6^\pm (q) & \la x \ra A_5^\pm (q) & (\la x \ra A_1(q)^3T_2).S_3 & q^{24}/6\\
E_7(q) & \la x \ra D_6(q) & (\la x \ra A_1(q)^6).S_3 & q^{48}/6 \\
E_8(q) & \la x \ra E_7(q) & (\la x \ra D_4(q)A_1(q)^3).S_3 & q^{96}/6 \\
\hline
\end{array}
\]

The number of pairs $(x^g,E)$ with $x^g \in E \in \D$ 
is of the order of $|\D|\cdot |D:(q+1)A_1(q)^3| \sim |\D|q^{18}$. Given $E$,
the proportion of conjugates $x^g \in E$ such that $\la x,x^g \ra = E$
is at least $1-c/q$ by Proposition \ref{d4gen}. 
Clearly $E$ is the unique member of $\D$ containing such 
$x^g$. Hence the number of conjugates $x^g$ such that $\la x,x^g\ra$
is a member of $\D$ is at least
$(1-c/q)|\D|q^{18}$. The number of conjugates of $x$ in
$G$ is $|G:C_G(x)|$, where $C_G(x)$ is as in the above table. Hence
the probability that $\la x,x^g\ra$ is a member of $\Delta$ is at least
\[
\frac{(1-\frac{c}{q})|\D|q^{18}}{|G:C_G(x)|} \ge \frac{1}{6} - \frac{c'}{q}.
\]
This completes the proof for $\e = +$. The proof for
$\e = -$ is similar. \hal 

\begin{prop} \label{g2probprop}
Let $G = G_2(q)$ with $q>2$ even; let $A_1$ and $\tilde A_1$ denote
commuting $SL_2(q)$ subgroups of $G$ generated by long 
and short root groups respectively;  
let $x,y$ be elements of order $3$ in $\tilde A_1, A_1$ respectively; 
and let $t$ be an involution in $A_1$.

\begin{enumerate}
\item[{\rm (i)}] $C_G(x) \cong SL_3^\e(q)$ and 
$C_G(y) = (q-\e)\times \tilde A_1$, 
where $q \equiv \e \bmod 3$.

\item[{\rm (ii)}] For random $g \in G$, the probability that 
$\la x, x^g \ra$ is a conjugate of 
$\tilde A_1$ is positive, and is at least $1-c_1/q$ where 
$c_1$ is an absolute constant.

\item[{\rm (iii)}] For random $g \in G$, the probability that 
$\la y, y^g \ra$ is a conjugate of $A_1$ is $\frac{1}{q^4(q^4+q^2+1)}$.

\item[{\rm (iv)}] For random $g \in G$, the probability that 
$\la A_1, t^g \ra \cong SL_3(q)$ is positive, and is at least $1/2-c_2/q$ 
where $c_2$ is an absolute constant.
\end{enumerate}
\end{prop}

\pf (i) If $L$ denotes the Lie algebra of type $G_2$,
then
\[
L \downarrow A_1\tilde A_1 = L(A_1\tilde A_1) \oplus (V(1)\otimes V(3)),
\]
where $V(i)$ denotes the irreducible module of
high weight $i$ (see, for example, \cite[11.12(ii)]{LSbk}). It follows
that $C_L(x)$, $C_L(y)$ have dimensions $8$ and $6$ respectively, so the 
centralizers of $x$ and $y$ in the
algebraic group $G_2$ are connected reductive subgroups of 
types $A_2$ and $T_1 A_1$.

(ii) This is similar to the proof of Proposition
\ref{probp2}.  Define
$\D = \{\tilde A_1^g : g \in C_G(x)\}$. Then
\[
|\D| = |C_G(x):C_G(x)\cap N_G(\tilde A_1)| = 
|SL_3^\e (q): A_1\cdot (q-\e)| \sim q^4.
\]
The number of pairs $(x^g,E)$ with $x^g \in E \in \D$ 
is of the order of $|\D|\cdot |\tilde A_1:(q-\e)| \sim q^{6}$. Arguing as in
Proposition \ref{d4gen}, we see that, given $E$, the proportion of
conjugates $x^g \in E$ such that $\la x,x^g \ra = E$ is at least
$1-c/q$. Clearly $E$ is the unique member of $\D$
containing such $x^g$. Hence the number of conjugates $x^g$ such that $\la
x,x^g\ra$ is a member of $\D$ is at least
$(1-c/q)q^6$, so the probability that $\la x,x^g\ra$ is
a member of $\D$ is at least
\[
\frac{(1-\frac{c}{q})q^6}{|G:C_G(x)|} \ge 1 - \frac{c'}{q}.
\]

The positivity statement in (ii) follows from the fact that 
$\tilde A_1 \cong SL_2(q)$ can be generated 
by two conjugates of $x$, which can be proved using (\ref{probbd}). 
Indeed, let $P$ be the probability that $\la x,x^g \ra = \tilde A_1$ for 
random $g \in \tilde A_1$. 
The maximal subgroups of
$\tilde A_1$ containing $x$ are 
$N:= N_{\tilde A_1}(\la x \ra) \cong D_{2(q-\e)}$ ($\e = \pm 1$) and 
conjugates of $SL_2(q_0)$ for maximal subfields $\F_{q_0}$ of $\F_q$. 
Observe that 
${\rm fix}_{\tilde A_1/N}(x) = 1$. 
If $M = SL_2(q_0)$ then $x^{\tilde A_1} \cap M = x^M$, so 
\[
{\rm fix}_{\tilde A_1/M}(x) = 
\frac{|\tilde A_1:M|\,|x^{\tilde A_1}\cap M|}{|x^{\tilde A_1}|} 
= \frac{|x^M|}{|x^{\tilde A_1}|} = \frac{q-\e}{q_0-\e_0},
\]
where $\e_0 = \pm1$ is such that $q_0 \equiv \e_0 \hbox{ mod }3$. Hence by (\ref{probbd}) 
\[
1-P \le \frac{2}{q(q+\e)} + \sum_{q_0} {\left(\frac{q-\e}{q_0-\e_0}\right)}^2 
\frac{q_0(q_0^2-1)}{q(q^2-1)},
\]
which is less than 1.

(iii) If $A_1^h$ is a conjugate of $A_1$ containing $y$, 
then $y,y^{h^{-1}} \in A_1$, so there exists $a\in A_1$
such that $y^{h^{-1}} = y^a$. Then $ah \in C_G(y)$ which is contained 
in $A_1\tilde A_1$ by (i), so 
$A_1^h = A_1$. In other words, the only conjugate of $A_1$ containing $y$ 
is $A_1$. Thus the 
probability in (ii) is $|y^G \cap A_1|/|y^G|$. The conclusion follows.

(iv) This is similar to the proof of (ii). Let $S$ be a subsystem subgroup
$SL_3(q)$ containing $A_1$, and let $\D = \{S^g:g \in \tilde
A_1\}$. Then $|\D| = |\tilde A_1: \tilde A_1 \cap N_G(S)|$ which is
at most $|\tilde A_1:(q+1).2|\sim q^2 / 2$. The
number of pairs $(t^g,E)$ with $t^g \in E \in \D$ is of the order of
$|\D|q^4 \sim q^6/2$, and a proportion of at least $1-c/q$ of
these satisfy $\la A_1, t^g \ra = E$. Since $|t^G| \sim q^6$, 
the lower bound in the 
conclusion follows. The positivity statement follows from the next proposition. \hal

\begin{prop}\label{3d4probprop}
Let $G = \,^3\!D_4(q)$ with $q$ even, and let $A$ be a long $SL_2(q)$
subgroup of $G$. Let $x$ and $t$ be elements of order $q+1$ and $2$ in $A$,
respectively. For random $g \in G$, the probability that $\la
x,t^g\ra$ is a subsystem subgroup $SL_3(q)$ is positive, and is at least
$1-c/q$, where $c$ is an absolute constant.
\end{prop}

\pf Let $S$ be a subsystem $SL_3(q)$ containing $A$, and let $\D =
\{S^g:g \in C_G(x) \}$. Note that $C_G(x) = \la x \ra SL_2(q^3)$ and
$|\D| = |SL_2(q^3):(q^3+1)| \sim q^6$. The number of pairs $(t^g,E)$
with $t^g \in E \in \D$ is of the order of $q^{10}$, and 
also $|t^G| \sim q^{10}$.
The lower bound follows as in the previous propositions.

For the positivity statement, let $S$ be a subsystem subgroup 
$SL_3(q)$ containing $x$. 
The maximal subgroups of $S$ appear in \cite[Tables 8.3--8.4]{BHRD}, 
from which we deduce that $x$ lies in 
just two maximal subgroups $P_1$, $P_2$, stabilizers of 1- and 2-spaces, 
respectively. 
These have structure
$(\F_q^2).(SL_2(q)\times (q-1))$, and each contains $q^3-1$ involutions. 
Since the total number of involutions in $S$
is $(q^3-1)(q+1)$, there is an involution $t$ such that 
$S = \la x,t \ra$, as required.
\hal

\section{Basic $SL_2$ subgroups in $SL_3(q)$ and 
$SL_6(q)$, $q$ odd}\label{sl3sl6}

Recall the definition of basic $SL_2$ subgroups in Section \ref{notnlie}. 
As components for our subsequent work, we require algorithms to
construct two basic $SL_2$ subgroups in a given
$SL_3(q)/Z$, and five basic
$SL_2$ subgroups in a given $SL_6(q)/Z$; here $q$ is odd and $Z$ is a
central subgroup.

In these and subsequent algorithms, we assume that our input group $G$ is 
described by a collection of generators in $GL_d(F)$ for some field 
$F$ of the same characteristic as $\F_q$.

\subsection{Algorithm for $SL_3(q)$}\label{sl3q}
Let $G$ be isomorphic 
to $SL_3(q)/Z$ with $q$ odd and $Z$ a central subgroup. The algorithm 
to construct two basic $SL_2$ subgroups in $G$ is the following.

\begin{enumerate}
\item[1.]
Find an involution $t_1 \in G$ by random search. 

\item[2.]
Construct $C_G(t_1)$ and $K_1 = C_G(t_1)' \cong SL_2(q)$. 

\item[3.] Find an involution $t_2 \in C_G(t_1)$ which does not
 commute with $K_1$, and compute $K_2 = C_G(t_2)' \cong SL_2(q)$.
\end{enumerate}
Now $K_1$ and $K_2$ are the required basic $SL_2$ subgroups of $G$. 

\begin{lemma} 
The algorithm is Monte Carlo, has probability greater than a positive absolute 
constant (independent of $q$) of finding the required involutions, and 
runs in polynomial time. 
\end{lemma}

\pf 
That we can both construct the involution $t_1$ and its
centralizer with positive probability independent of $q$ 
follows from Section
\ref{centralizer}. Now consider the second involution $t_2$.
There is a maximal torus $T$ of order $(q-1)^2/|Z|$ in
$C:=C_G(t_1)$, and at least 1/4 of its regular elements power
into the conjugacy class in $C$ of a suitable involution
$t_2$. The number of non-regular elements in $T$ is at most $3(q-1)$,
and $|N_C(T):T| = 2$. By Lemma \ref{parwil}, the proportion of
elements of $C$ which power to a suitable involution $t_2$ is at least
$\frac{1}{2}(\frac{1}{4} - \frac{3}{q-1})$. Therefore this is a lower bound
for the probability of finding $t_2$ and is 
positive for every $q$ since $t_2$ exists.
In Section \ref{prelim} we cite polynomial-time
algorithms to perform the other tasks.
 \hal

\subsection{Algorithm for $SL_6(q)$ }\label{sl6q}

Let $G$ be isomorphic to
$SL_6(q)/Z$, where $q$ is odd and $Z$ is a central subgroup. 
Involutions in $G$ have
centralizers with derived groups $SL_4(q)\circ SL_2(q)$, $SL_5(q)$,
$SL_3(q)\circ SL_3(q)$ or $SL_3(q^2)/Z$. In our algorithm we 
consider only involutions having centralizers of the first type, and we
call such a centralizer ``good"; we can inspect orders  of random
elements in the involution centralizer to determine whether 
the centralizer is good. 

The algorithm 
to construct five basic $SL_2$ subgroups in $G$ is the following.

\begin{enumerate}
\item[1.] Find an involution $t_1 \in G$ with good centralizer,
so $C_G(t_1)' = SL_4(q)\circ SL_2(q)$. Use \KillFactor\
(see Section \ref{prelim}) to construct the factor $K_1 \cong SL_2(q)$ of this centralizer.

\item[2.]
Find an involution $t_2 \in C_G(t_1)$ with good centralizer
such that $[t_2,K_1]\ne 1$. Construct $K_2$, 
the $SL_2(q)$ factor of $C_G(t_2)'$. 

\item[3.] Find an involution $t_3 \in C_G(t_1,t_2)$ with good
 centralizer such that $t_3 t \not \in Z(G)$ for all $t \in \la
 t_1,t_2\ra$, and $[t_3,K_1]=1$, $[t_3,K_2] \ne 1$. Construct $K_3 =
 C_G(t_1,t_2,t_3)' \cong SL_2(q)$.

\item[4.] Find an involution $t_4 \in C_G(t_1,t_2,t_3)$ with good
 centralizer such that $t_4 t \not \in Z(G)$ for all $t \in \la
 t_1,t_2,t_3\ra$, and $[t_4,K_1]=[t_4,K_2]=1$, $[t_4,K_3] \ne
 1$. Construct $K_4$, the $SL_2(q)$ factor in $C_G(t_4)'$.

\item[5.]
Let $t_5 = t_1t_3$ and construct $K_5$, the $SL_2(q)$ factor in $C_G(t_5)'$. 
\end{enumerate}
With respect to a suitable basis of $V_6(q)$, 
$\pm t_1 = (-1,-1,1,1,1,1)$, $\pm t_2 = (1,-1,-1,1,1,1)$,
$\pm t_3 = (1,1,-1,-1,1,1)$, and $\pm t_4 = (1,1,1,-1,-1,1)$. 
Hence $K_1,\ldots,K_5$ are the required basic $SL_2$ subgroups. 

\begin{lemma} 
The algorithm is Monte Carlo, has probability greater than a positive absolute constant 
of finding the required involutions, and 
runs in polynomial time. 
\end{lemma}
\pf 
In Steps 2, 3 and 4 we use a maximal torus of
order $(q-1)^5$ and Lemma \ref{parwil} to estimate the
probabilities of finding suitable involutions $t_2,t_3,t_4$. We
illustrate the calculation for $t_2$. Write $t_1 = (-1,-1,1,1,1,1)$ as
above, and let $T$ consist of the diagonal matrices $(\a_1,\ldots
,\a_6)$ where $\a_i \in \F_q^*$ and $\prod \a_i = 1$. Let $Q$ be the subgroup
of index 2 in $\F_q^*$. If we take $\a_2 \in \F_q^* \setminus Q$,
$\a_1 \in Q$ and the other $\a_i$ arbitrary, then this element of $T$
powers to a suitable involution $t_2$, and the number of such elements
in $T$ is $|T|/4$, of which at most $f(q)$ are non-regular, for some
polynomial $f(q)$ of degree at most 4. Also $|N_C(T):T| = 48$, where
$C = C_G(t_1)$. Hence Lemma \ref{parwil} shows that the proportion of
elements of $C$ powering to a suitable involution $t_2$ is at least
$1/192 - c/q$ for some absolute constant $c$.
\hal

\section{Basic $SL_2$ subgroups in $E_6(q)$, $E_7(q)$ 
and $E_8(q)$, $q$ odd}\label{e678sec}

Let $G$ be isomorphic to one of
the quasisimple groups of type $G(q) = E_6(q)$, $E_7(q)$ 
or $E_8(q)$ with $q$ odd. 
We first present algorithms 
to construct basic $SL_2$ subgroups of $G$ 
and later justify them. Each algorithm starts with 
the construction of an involution 
centralizer; these are 
described in Proposition \ref{invcent}. 

As usual, we assume that our input group $G$ is 
described by a collection of generators in $GL_d(F)$ for some field 
$F$ of the same 
characteristic as $\F_q$.

\subsection{$E_6(q)$, $q$ odd}\label{e6}
Let $G(q) = E_6(q)$, $q$ odd. 

\begin{enumerate}
\item[1.] Find an involution $t_0 \in G$ with $C_G(t_0)$ of type $A_1A_5$. 
Construct the factors $K_0 \cong SL_2(q)$ and $D \cong A_5(q)$ of $C_G(t_0)'$.

\item[2.] Find an involution 
$t_2 \in C_G(t_0)$ such that $C_D(t_2)' = 
E \cong SL_3(q)\circ SL_3(q)$. Construct the two $SL_3(q)$ factors.

\item[3.] 
Construct basic $SL_2$ subgroups $K_1,K_3$
 in the first factor $SL_3(q)$ of $E$, and $K_5,K_6$ in 
the second factor. Let $Z(K_i) = \la t_i \ra$.

\item[4.] Construct $K_4$, the $SL_2(q)$ factor 
of $C_D(t_1,t_6)' \cong SL_2(q)^3$ which centralizes $K_1K_6$. 
Now $K_1,K_3,K_4,K_5,K_6$ are basic 
$SL_2$ subgroups in $D \cong SL_6(q)$.

\item[5.] The centralizer $C_G(t_2)$ is of type $A_1A_5$; 
construct $K_2$, the $SL_2(q)$ factor of $C_G(t_2)$. 
\end{enumerate}

We now have the six basic $SL_2$ subgroups 
$K_1,\ldots, K_6$ in the $E_6$ Dynkin diagram: 
\[
\begin{array}{llllcllll}
K_1& - &K_3& - &K_4& - & K_5& - &K_6 \\
 && && | & & & & \\
 &&& & K_2 & &&& 
\end{array}
\]

\subsection{$E_7(q)$, $q$ odd} \label{e7}
This algorithm is similar to that for $E_6$. 

\begin{enumerate}
\item[1.] 
Find an involution $t_0 \in G$ such that 
$C_G(t_0)$ is of type $A_1D_6$, and construct the factors $K_0 \cong SL_2(q)$ 
and $D \cong D_6(q)$ of $C_G(t_0)'$.

\item[2.] 
Find an involution $t_1 \in C_G(t_0)$ such that $C_D(t_1)' = E \cong A_5(q)$
and $C_G(t_1)$ is of type $A_1D_6$.

\item[3.]
Construct basic $SL_2$ subgroups 
$K_2$, $K_4$, $K_5$, $K_6$, $K_7$ in $E$. 
Let $t_i$ be the central involution in $K_i$.

\item[4.]
The element $t_1$ is a root involution; construct the factor 
$K_1 \cong SL_2(q)$ of $C_G(t_1)$. 

\item[5.]
The element $t_3 = t_0t_5t_7$ is a root involution; 
construct the factor $K_3\cong SL_2(q)$ of $C_G(t_3)$.
\end{enumerate}

We now have the seven basic $SL_2$ subgroups 
$K_1,\ldots, K_7$ in the $E_7$ Dynkin diagram: 
\[
\begin{array}{llllcllllll}
K_1& - &K_3& - &K_4& - & K_5& - &K_6& - & K_7\\
 &&& & | &&&& & &\\
 & && & K_2 & & &&&&
\end{array}
\]

\subsection{$E_8(q)$, $q$ odd}\label{e8}
This algorithm is similar to that for $E_6$ and $E_7$. 

\begin{enumerate}
\item[1.] 
Find an involution $t_0 \in G$ with $C_G(t_0)$ of type $A_1E_7$, 
and construct the factors 
$K_0 \cong SL_2(q)$, $D \cong E_7(q)$ of $C_G(t_0)'$.

\item[2.]
Find an involution $t_8 \in C_G(t_0)$ such that $C_D(t_8)' = E \cong E_6(q)$.

\item[3.] 
Construct basic $SL_2$ subgroups $K_1,\ldots K_6$ in $E$. 

\item[4.]
The element $t_8$ is a root involution; construct the 
factor $K_8 \cong SL_2(q)$ of $C_G(t_8)$.

\item[5.]
The element $t_7 = t_0t_2t_5$ is a root involution; 
construct the factor $K_7\cong SL_2(q)$ of $C_G(t_7)$.
\end{enumerate}

We now have the eight basic $SL_2$ subgroups 
$K_1,\ldots, K_8$ in the $E_8$ Dynkin diagram: 
\[
\begin{array}{llllcllllllll}
K_1& - &K_3& - &K_4& - & K_5& - &K_6& - & K_7& - &K_8\\
 &&& & | &&&&& & & &\\
 &&& & K_2 & & & &&&&&
\end{array}
\]

\subsection{Justification}\label{juste678}

\begin{prop}\label{e678}
The algorithms for $E_6(q)$, $E_7(q)$ and $E_8(q)$ for odd $q$ described
above are Monte Carlo and run in polynomial time.
\end{prop}

\pf We first prove the correctness of the algorithm for $E_6(q)$.
In Step 1, finding the involution $t_0$ and constructing its centralizer
is justified by the results of Section \ref{centralizer}.
The factors $K_0$ and $D$ of $C_G(t_0)'$ are
constructed using the algorithm \KillFactor, referred to in Section
\ref{prelim}.

Now consider Step 2 of the algorithm: find an involution $t_2 \in
C_G(t_0)$ with centralizer containing $SL_3(q) \circ SL_3(q)$. We 
show that there is a positive lower bound (independent of $q$) for 
the probability of
finding such an involution. This does not follow directly from
Theorem \ref{powerup}, but follows from the method of its 
proof in \cite{PW}. Namely, there is a maximal torus of $C_G(t_0) \cong
(SL_2(q)\circ A_5(q)).2$ of order $(q^3-1)^2$, and at least 1/8
of the elements of this torus power to involutions which have the
desired centralizer structure; thus Lemma \ref{parwil} gives the
required conclusion. The centralizer can be computed by 
Section \ref{centralizer},
and the $SL_3(q)$ factors extracted using \KillFactor.

Step 3 of the algorithm is justified in Section \ref{sl3sl6}. In
Step 4, the construction of $C_D(t_1)$, and of the centralizer of
$t_6$ within this group, is justified using Section \ref{centralizer}. 
Observe that if $\la t_4 \ra = Z(K_4)$, then $t_4 = t_1t_6t_0$.

Step 5 requires a little more argument.
Recall that $E \cong SL_3(q) \circ SL_3(q)$,
a central product of two $SL_3(q)$ subsystem subgroups of $G$. From
the subsystem $A_2^3$ of the $E_6$ root system, we see that $C_G(E)$
is isomorphic to $Z(E)SL_3(q)$, so $t_2$ and $K_0$ are contained
in $C: = C_G(E)' \cong SL_3(q)$. Hence $t_2$ is a root involution and 
we let $K_2$ be the $SL_2(q)$ factor of $C_G(t_2)$.

Finally, we show that $K_1,\ldots, K_6$ pairwise generate either
their direct product or $SL_3(q)$ according to their positions in the
Dynkin diagram. Observe that $C_G(t_2)' = K_2S$ with $S \cong SL_6(q)$, and
clearly $E < S$. Hence $K_2 \le C_G(E)' = C$. Therefore $K_2$
centralizes each of $K_1,K_3,K_5,K_6$ and $\la K_2,K_0 \ra = C$. Hence 
$K_0,K_2$ and $K_4$ are contained in $C_G(K_1K_6)' \cong SL_4(q)$. The
central involutions $t_0,t_2,t_4$ commute: $t_2$ commutes with $t_4$
since it commutes with $t_0,t_1,t_6$ and $t_4 = t_1t_6t_0$. Working
in $SL_4(q)$ relative to a basis diagonalizing these three
involutions, we see that $\la K_2,K_4\ra \cong SL_3(q)$.
This justifies the algorithm for $E_6(q)$. 

For $E_7(q)$ the proof is similar, with the following additional
observations. In Step 2, such an involution $t_1$ can be found with
positive probability by the usual argument using Lemma \ref{parwil} 
and a maximal torus of order $(q- 1)(q^6-1)$ in $C_G(t_0)$. In Step
4, observe that $K_1$ commutes with $E$ and $\la K_0,K_1 \ra = C_G(E)
\cong SL_3(q)$. For Step 5 and the Dynkin diagram generation of
the $K_i$, observe the equation between toral elements $h_0(-1) =
h_{2234321}(-1) = h_3(-1)h_5(-1)h_7(-1)$ (using the 
notation of Section \ref{notnlie}).
Since $t_i = h_i(-1)$, the involution in the centre of the
final $SL_2(q)$ to complete the Dynkin diagram must be $t_3 =
t_0t_5t_7$. 

Finally consider $E_8(q)$.
To justify Step 2 we take a maximal torus of order $(q\pm 1)^2(q^6+q^3+1)$ 
in $C_G(t_0)$ and use Lemma \ref{parwil} as usual. In Step 5
the equation $h_0(-1) = h_{23465432}(-1) = h_2(-1)h_5(-1)h_7(-1)$
justifies our definition $t_7 = t_0t_2t_5$, and implies that the
factor $K_7\cong SL_2(q)$ of its centralizer completes the Dynkin
diagram as claimed. \hal

\section{Basic $SL_2$ subgroups in $E_6(q)$, $E_7(q)$ and $E_8(q)$, 
$q$ even}\label{e678p2sec}

In this section we assume that $G$ is described by a collection of generators
in $GL_d(F)$, and $G$ is quasisimple and isomorphic to one of $G(q) =
E_6(q)$, $E_7(q)$ or $E_8(q)$, where $F$ and $\F_q$ are both finite
fields of characteristic 2, and $q>2$. 
We first present algorithms 
to construct basic $SL_2$ subgroups of $G$,
and later justify them. 

Throughout, $\o$ denotes a generator for the multiplicative group of $\F_q$.

\subsection{$E_6(q)$, $q$ even} \label{e6p2sec}
Assume $G$ is isomorphic to $E_6(q)$ with $q$ even and $q>2$. 
\begin{enumerate}
\item[1.]
Find $y \in G$ of order $(q+1)(q^5-1)/d$ where $d=(3,q-1)$, 
and define $x = y^{(q^5-1)/d}$. 

\item[2.] Find $g \in G$ such that $X := \la x,x^g\ra$ is isomorphic to
 $D_4^\e(q)$ (where $\e = \pm$). 

\item[3.] Construct an isomorphism $\phi$ from $X$ to the standard copy of
 $D_4^\e (q) = \O(V)$, where $V = V_8(q)$.

\item[4.] Find a standard basis $e_1,e_2,e_3,f_3,f_2,f_1$ for 
a non-degenerate subspace of $V$ of type $O_6^+$.
In the $SL_3(q)$ subgroup of $X\phi$ fixing $\la
 e_1,e_2,e_3\ra$ and $\la f_1,f_2,f_3\ra$, write down six elements acting
 on $\la e_1,e_2,e_3\ra$ as $v_i,u_i^+,u_i^-$ ($i=1,2$), where
\[
\begin{array}{l}
v_1 = \pmatrix{\o^{-1} & 0&0 \cr 0&\o & 0 \cr 0&0&1},\;
u_1^+ = \pmatrix{1 & 1&0 \cr 0&1 & 0 \cr 0&0&1},\;
u_1^- = \pmatrix{1 & 0&0 \cr 1&1 & 0 \cr 0&0&1}, \\
v_2 = \pmatrix{1& 0&0 \cr 0&\o^{-1} & 0 \cr 0&0&\o },\;
u_2^+ = \pmatrix{1 & 0&0 \cr 0&1 & 1 \cr 0&0&1},\;
u_2^- = \pmatrix{1 & 0&0 \cr 0&1 & 0 \cr 0&1&1}.
\end{array}
\]
Abusing notation, write also $v_i,u_i^+,u_i^-$ for the
inverse images of these elements under $\phi$. Define $K_0 = \la
v_1,u_1^+,u_1^- \ra$ and $K_2 = \la v_2,u_2^+,u_2^- \ra$, basic
$SL_2$ subgroups in $X$.

\item[5.] 
Construct the involution centralizer $C_G(u_1^+)$. 

\item[6.] Apply Lemma \ref{Formula} and the ensuing remark to 
$N: = \la C_G(u_1^+),v_1\ra = Q(D\times \la v_1 \ra)$, where $D =
 C_G(K_0) \cong A_5(q)$ and $Q = O_2(N)$. 
This constructs $C_N(v_1) = (Z(Q) \times D)\la v_1 \ra$. 
Construct its second derived group, $D$.

\item[7.] Construct an isomorphism $\psi$ from $D$ to the standard copy of
 $SL_6(q)= SL(V)$ modulo a central subgroup $Z$ of order either 1 or
 $(3,q-1)$. (Here $|Z|=3$ if and only if $3|q-1$ and $Z(G)=1$.) 

\item[8.] Consider $v_2 \in K_2$. This element acts on $D$. 
Compute $T \in GL_6(q)$ such
 that $(g^{v_2})\psi = (g\psi)^T$ for all $g \in D$.  

\item[9.] Diagonalise $T$ to find a basis of $V$ with respect to
 which $T = (\o I_3,\o^{-1}I_3)$. Let this basis be $x_1,\ldots x_6$.

\item[10.] For $1\le i\le 5$, let $a_i,b_i,c_i$ be the inverse images
 under $\psi$ of matrices fixing $x_j$ for $j\ne i,i+1$ and such that
\[
\begin{array}{l}
a_i \psi \,:\,x_i \rightarrow\o^{-1} x_i,\, x_{i+1} \rightarrow \o x_{i+1}, \\
b_i \psi \,:\, x_i \rightarrow x_i+x_{i+1}, \,x_{i+1}\rightarrow x_{i+1}, \\
c_i \psi \,:\, x_i \rightarrow x_i, \,x_{i+1}\rightarrow x_i+x_{i+1}.
\end{array}
\]
Define $K_1 = \la a_1,b_1,c_1\ra$, and for $i=3,\ldots, 6$, define
$K_i = \la a_{i-1},b_{i-1},c_{i-1}\ra$. Then $K_1,K_3,K_4,K_5,K_6$
are basic $SL_2$ subgroups in $D \cong SL_6(q)$.
\end{enumerate}

We now have the six basic $SL_2$ subgroups $K_1,\ldots, K_6$ in
the $E_6$ Dynkin diagram.

\subsection{$E_7(q)$, $q$ even} \label{e7p2sec}

Assume $G$ is isomorphic to $E_7(q)$ with $q$ even and $q>2$. 
\begin{enumerate}
\item[1-6.]
These steps are as for the $E_6$ algorithm, with the following modifications. 
In Step 1, we find an element $y$ of order $(q+1)(q^5-1)$ and 
define $x = y^{q^5 - 1}$; 
the basic $SL_2$ subgroups constructed in Step 4 are $K_0$ and $K_1$;
in Step 6, we construct $D = C_G(K_0)$ which is isomorphic to $D_6(q)$.

\item[7.] Construct an isomorphism $\psi$ from $D$ to the standard copy of
 $D_6(q) = \O_{12}^+(q)= \O (V)$, where $V$ has associated bilinear
 form $(\,,\,)$. 

\item[8.] Consider $v_2 \in K_1$. This element acts
 on $D$. Compute $T \in GL_{12}(q)$ 
such that $(g^{v_2})\psi = (g\psi)^T$ for all $g \in D$.

\item[9.] Diagonalise $T$ to find a basis of $V$ with respect to
 which $T = (\o I_6, \o^{-1}I_6)$. Choose a basis $e_1,\ldots, e_6$
 for the $\o$-eigenspace, and a basis $f_1,\ldots, f_6$ for the
 $\o^{-1}$-eigenspace such that $(e_i,f_j) = \d_{ij}$.

\item[10.] 
 For $1\le i\le 5$, let $a_i,b_i,c_i \in D$ be the inverse 
images under $\psi$ of
 the matrices in $\O (V)$ fixing $e_j,f_j$ for $j\ne i,i+1$ and such that
\[
\begin{array}{l}
a_i \psi \,:\,e_i \rightarrow\o^{-1} e_i, e_{i+1} \rightarrow \o e_{i+1}, 
f_i \rightarrow\o f_i, f_{i+1} \rightarrow \o^{-1} f_{i+1}, \\
b_i \psi \,:\, e_i \rightarrow e_i+e_{i+1}, e_{i+1}\rightarrow e_{i+1}, 
f_i \rightarrow f_i, f_{i+1}\rightarrow f_{i+1}+f_i, \\
c_i \psi \,:\, e_i \rightarrow e_i, e_{i+1}\rightarrow e_i+e_{i+1}, 
f_i \rightarrow f_i+f_{i+1}, f_{i+1}\rightarrow f_{i+1}.
\end{array}
\]
Define $K_2 = \la a_1,b_1,c_1\ra$, and for $i=4,\ldots, 7$, define
$K_i = \la a_{i-2},b_{i-2},c_{i-2}\ra$. Finally, define $K_3 = \la
a_6,b_6,c_6\ra$, where these elements are the inverse images of the
matrices fixing $e_j,f_j$ for $j \ge 3$ and such that
\[
\begin{array}{l}
a_6 \psi \,: \, e_1 \rightarrow \o e_1, e_2 \rightarrow \o e_2, 
f_1 \rightarrow \o^{-1}f_1, f_2 \rightarrow \o^{-1}f_2, \\
b_6 \psi \,:\, e_1\rightarrow e_1+f_2, e_2\rightarrow e_2+f_1, 
f_1\rightarrow f_1,f_2\rightarrow f_2, \\
c_6 \psi \,:\, e_1\rightarrow e_1, e_2\rightarrow e_2, f_1\rightarrow f_1+e_2,
f_2\rightarrow f_2+e_1.
\end{array}
\]
Then $K_2,K_3,K_4,K_5,K_6,K_7$ are basic $SL_2$ subgroups
 in $D \cong D_6(q)$.
\end{enumerate}

We now have the seven basic $SL_2$ subgroups $K_1,\ldots, K_7$
in the $E_7$  Dynkin diagram.

\subsection{$E_8(q)$, $q$ even} \label{e8p2sec}

Assume $G$ is isomorphic to $E_8(q)$ with $q$ even and $q>2$. 
\begin{enumerate}
\item[1-6.] These steps are 
as for the $E_6$ algorithm, with
 the following modifications. In Step 1, we find an element $y$ of order
 $(q+1)(q^7-1)$ and define $x = y^{q^7-1}$; 
the basic $SL_2$ subgroups constructed in Step 4 are $K_0$ and $K_8$;
in Step 6, we construct $D = C_G(K_0)$ which is isomorphic to $E_7(q)$.

\item[7.] 
Using the algorithm of Section \ref{e7p2sec},
construct basic $SL_2$ subgroups $K_1,\ldots, K_7$ of $D$; 
 label root elements of $D$ as in Section \ref{e678labsec}. 
Construct an isomorphism $\psi$ from $D$ to the 
standard copy of $E_7(q)$ (a group of $56\times 56$ matrices).

\item[8.] Consider $v_2 \in K_8$. 
This element acts  on $D$. 
Compute $T \in  D\psi$ such that 
$(g^{v_2})\psi = (g\psi)^T$ for all $g \in D$.

\item[9.] 
Compute $g \in D\psi$ such that $T^g =  h_{2346543}(\l)$ for 
some $\l \in \F_q$. 
Replace $K_8$ by $K_8^{g\psi^{-1}}$.
\end{enumerate}

We now have the eight basic 
$SL_2$ subgroups $K_1,\ldots, K_8$ in the $E_8$ Dynkin diagram.

\subsection{Justification}

\begin{prop}\label{e678propp2} 
The algorithms for $E_6(q)$, $E_7(q)$ and $E_8(q)$  for even $q$ 
described above are Monte Carlo and run in polynomial time.
\end{prop}

\pf 
In Step 1 of each algorithm, the justification for being able to 
find an element $y$ of the specified order 
is standard.  Consider, for example, the $E_6$ case: 
in the simple group $G := E_6(q)$ there is a cyclic maximal torus $T$ of order
$t:=(q+1)(q^5-1)/d$ (see \cite[\S 2]{KS}); the number of generators of 
$T$ is $\phi(t) > t/c\log \log q$ where $c$ is an absolute constant; 
hence the proportion of elements of order $t$ in $G$ is at least 
$1/(c\log \log q\cdot |N_G(T):T|)$.

Observe that $y$ lies in a
maximal torus of $G$ contained in a subsystem subgroup of type
$A_1A_5$, $A_1D_6$ or $A_1E_7$ (see \cite[\S 2]{KS}).
Hence the power $x$ of $y$ must
lie in a long $SL_2(q)$ subgroup of $G$. By Proposition
\ref{probp2}, there is a 
positive lower bound independent of $q$ for the probability 
that, for random $g \in G$, 
$\la x,x^g\ra$ is $D_4^\e (q)$, a subsystem subgroup, as required for Step 2.

In Step 3, the construction of an isomorphism 
$X \rightarrow D_4^\e(q)$ is justified by
Theorem \ref{available-crec}.

In Step 5, the element $u_1^+$ is a root involution so the construction of 
the involution centralizer $C_G(u_1^+)$
is justified by Theorem \ref{even-cent} and the ensuing remark.
The structure of $C_G(u_1^+)$ is given in \cite{AS}:
$C_G(u_1^+) = QD$, where $D = C_G(K_0) \cong A_5(q)$,
$D_6(q)$ or $E_7(q)$, and
$Q \cong q^{1+20}$, $q^{1+32}$ or $q^{1+56}$,
when $G$ is of type $E_6,E_7,E_8$ respectively.
The element $v_1$
normalizes $C_G(u_1^+)$ (which is $C_G(U)$ where $U$ is a
root group of $K_0$ containing $u_1^+$), centralizes $D$, and acts
fixed point freely on $Q/Z(Q)$. Lemma \ref{Formula} and the ensuing 
remark gives a construction of
$(Z(Q) \times D)\la v_1 \ra$, as claimed in Step 6.

For $E_6(q)$ (or $E_7(q)$), the isomorphism in Step 7 from $D$ to
$A_5(q)$ (or $D_6(q)$) is justified by 
Theorem \ref{available-crec}.
The remaining steps ensure that 
$v_2$ acts as $h_2(\o)$ (or $h_1(\o)$) on $D$. Hence we choose
the remaining $SL_2(q)$ subgroups in $D$ to fit in with the subgroups
$K_0,K_2$ (or $K_0,K_1$) already defined. 
That they pairwise generate the correct
groups is established as in the proof of Proposition \ref{e678}.  
In Step 8, the computation of
the matrix $T$ involves solving linear equations in the entries of $T$
of the form $TA = BT$, where $A = (g^{v_2})\psi$, $B = g\psi$ for
generators $g$ of $D$; such systems of equations
over $\F_q$ can be solved in polynomial time.

For the $E_6$ case, in $SL_3(q) \cong \la K_2,K_4\ra$, $K_2$
and $K_4$ satisfy the correct picture of being the subgroups $(X,1)$
and $(1,X)$, for $X \in SL_2(q)$, relative to some basis of the natural
module $V = V_3(q)$: indeed, we constructed $K_4$ so that it is stabilized by
$v_2 \in K_2$, which implies that $v_2$ stabilizes $C_V(K_4)$; thus
$C_V(K_4) \subseteq [V,K_2]$ as required. 
Similar remarks apply in the $E_7$ case to 
$SL_3(q) \cong \la K_1,K_3\ra$.

For $E_8(q)$, Steps 7--9 are more complex. 
In Step 7, we construct an isomorphism $\psi$ from $D$ to the standard copy of
$E_7(q)$, a specific group of $56 \times 56$ matrices with standard 
generators $\hat {\cal S}$. 
Specifically, we find, 
as in Sections \ref{e678pres} 
a set ${\cal S}$ of standard generators 
of $D$ which satisfies the reduced Curtis-Steinberg-Tits presentation 
of $E_7(q)$.
For $x \in D$, we use the algorithm of \cite{CMT}, 
applied to the action of $D$ on
an absolutely irreducible composition factor of $V_d(F)\downarrow D$, 
to express $x$ as a word $w({\cal S})$; 
now $\psi$ is defined to send $x$ to $w(\hat {\cal S})$. 
That this step can be performed in polynomial time follows from
this proposition (already proved for $E_7(q)$), together with   
the algorithms of Section \ref{e678pres} and
\cite{CMT}. 

Since $D = C_G(v_1)'$, the element $v_2$ acts on $D$, and induces an
inner automorphism. In Step 8, we use linear algebra to find 
$T' \in GL_{56}(q)$ such that $(g^{v_2})\psi = (g\psi)^{T'}$ for all 
generators $g$ of $D$. 
Some scalar multiple, $T$, of $T'$ of determinant 1 must 
lie in $D\psi$; we use 
\cite{CMT} to determine $T$.

The centralizer of $T$ in $D\psi$ contains the image under $\psi$ of 
$C_D(K_8) = C_G(K_0,K_8) \cong E_6(q)$. 
It follows that $T$ is $D\psi$-conjugate to the toral element 
$h:=h_{2346543}(\l)$ for some eigenvalue $\l$ of $T$ on $V_{56}(q)$.
We compute $g \in D\psi$ conjugating $T$ to $h$ as follows.

\begin{enumerate}
\item 
Map $D\psi$ to its action on the Lie algebra $L$ of type $E_7$ over $\F_q$. 
Call this map $\phi$.
\item 
In each of $C_L(T\phi)$ and $C_L(h\phi)$, 
compute a split Cartan subalgebra by taking the 
centralizer  of a random semisimple
element. We claim that this is a split Cartan subalgebra with 
probability at least 
$c(1-|\Phi|/q)$, where $\Phi :=\Phi(E_6)$, the $E_6$
root system, and $c$ is a positive absolute constant. 
To prove the claim, observe that
\[
C_L(h\phi) = \la z\ra \oplus L(E_6) = H \oplus \sum_{\a \in \Phi(E_6)} L_\a
\]
where $z$ is semisimple, $H$ is a split Cartan subalgebra of $L$, and 
the $L_\a$ are 1-dimensional 
root spaces for $\a \in\Phi$. If $v \in H$ satisfies 
$\a (v) \ne 0$ for all $\a$, then $C_{C_L(h\phi)}(v) = H$ -- that is, $v$ is 
regular semisimple in $C_L(h\phi)$. 
The number of such $v\in H$ is at least $|H|(1-|\Phi|/q)$, so the total number 
of regular semisimple
elements in $C_L(h\phi)$ is at least this number multiplied by the number of 
conjugates of $H$ under the group
$C_{D\psi \phi}(h\phi) \cong E_6(q) \circ (q-1)$. For large $q$, the stabilizer of $H$ in the 
latter group is the normalizer of a
Cartan torus, of order $(q-1)^7|W(E_6)|$. It follows that the number of 
regular semisimple elements in $C_L(h\phi)$
is at least
\[
(|E_6(q) : (q-1)^6 W(E_6)|) \cdot |H|(1-|\Phi|/q),
\]
which is at least $c(1-|\Phi|/q)\cdot |C_L(h\phi)|$. This proves the claim. 
That we have found a split Cartan subalgebra can 
be verified in polynomial time by the argument of \cite[5.2]{CM09}. 

\item 
Use the polynomial-time algorithm of \cite[Theorem 1]{CR} 
to compute Chevalley bases $B_T$, $B_h$ of $L$ with respect to the 
Cartan subalgebras constructed in Step 2.

\item The element $g'$ of $GL(L)$ 
conjugating $B_T$ to $B_h$ lies in $D\psi\phi$, and conjugates 
$T\phi$ to an element of a Cartan torus of $D\psi \phi$ containing $h\phi$.
\item 
Adjust $g'$ by a computation in the Weyl group of $E_7(q)$ 
to an element $g''$ of $D\psi\phi$ conjugating $T\phi$ to $h\phi$. 
Take $g = g''\phi^{-1} \in D\psi$, as required.
\end{enumerate}

For convenience, we now abuse notation and write $g$ instead of $g\psi^{-1}$.
To complete the proof, we argue that replacing $K_8$ by $K_8^g$ provides
a set $K_1,\ldots, K_8$ of basic $SL_2$ subgroups.
For this, we need only to check that $K_8^g$ 
centralizes $K_1,\ldots ,K_6$ and 
$\la K_7,K_8^g \ra \cong SL_3(q)$. 
First observe that $v_2^g \in N_G(D) = DK_0$, 
so $v_2^g = hk_0$ with $k_0 \in K_0$. Also $C_G(K_8^g) = 
C_G(v_2^g)' \ge C_D(h)$, and this contains 
$K_1,\ldots ,K_6$. Finally $C_G(K_7) = C_G(h_7)'$ where $h_7 = h_{\a_7}(\o)$, so 
$C_G(K_7,K_8^g) = C_G(h_7,hk_0)'$. 
We claim that this centralizer is of type $E_6(q)$.
Indeed, we can label the $E_8$ root system so that 
$k_0 = h_{\a_0}(\mu) = h_{23465432}(\mu)$
for some $\mu \in \F_q$; the fact 
that $hk_0$ is conjugate to $v_2$ forces $\mu = \l^{-1}$ or 
$\l^{-3}$ (recall that $h=h_{2346543}(\l)$). 
Now considering $h_7$ and $hk_0$ as elements of the
subsystem subgroup $A_3$ corresponding to the 
roots $\a_7,\a_8,\a_0$, we see that they lie in an
$A_2$ subsystem, and hence centralize an $E_6$ subsystem in $E_8$. 
This proves the claim.
Hence  $\la K_7,K_8^g \ra \le C_G(E_6(q)) \cong SL_3(q)$.
Since $\la K_7,K_8^g\rangle$ contains $\la h_7,hk_0 \ra$, 
a toral subgroup of rank 2, it follows that  
$\la K_7,K_8^g \ra \cong SL_3(q)$ as required. \hal

\section{Basic $SL_2$ subgroups in $F_4(q)$}\label{f4sec}

\subsection{$F_4(q)$, $q$ odd}
Let $G$ be isomorphic to $F_4(q)$ with $q$ odd.
We present an algorithm to construct basic $SL_2$ subgroups in $G$.

\begin{enumerate}
\item[1.] 
Find an involution $t_0 \in G$ such that 
$C_G(t_0) \cong (SL_2(q)\circ Sp_6(q)).2$. Construct the factors
 $K_0\cong SL_2(q)$ and $D\cong Sp_6(q)$ of the centralizer.

\item[2.]
Find an involution $t_1 \in C_G(t_0)$ such that $C := C_D(t_1)' \cong SL_3(q)$.

\item[3.] 
Construct basic $SL_2$
 subgroups $K_3,K_4$ in $C$. Let $t_i$ be the involution in $K_i$. 

\item[4.]
Let $t_2 = t_0t_4$, a root involution; 
construct $K_2$, the $SL_2(q)$ factor in $C_G(t_2)$.

\item[5.]
Also $t_1$ is a root involution; construct $K_1$, the $SL_2(q)$
factor of its centralizer.
\end{enumerate}

We now have the four basic $SL_2$ subgroups 
$K_1,K_2,K_3,K_4$ in the Dynkin diagram:
\[
\begin{array}{lllllll}
K_1 & - & K_2 & =>= & K_3& - & K_4
\end{array}
\]

\begin{prop}\label{f4prop}
The above algorithm for $F_4(q)$ for odd $q$ is Monte Carlo and 
runs in polynomial time.
\end{prop}

\pf Finding the involutions and centralizers 
in Steps 1 and 2 is justified in the usual way 
using Lemma \ref{parwil} and Section \ref{centralizer}.
For Step 4, with respect to a suitable basis for 
the natural 6-dimensional module for $D \cong Sp_6(q)$, 
$t_3 = (-1,-1,-1,-1,1,1)$, $t_4 = (-1,-1,1,1,-1,-1)$ and $t_0 = -I$; 
hence $t_2 = t_0t_4$ is a root involution.

Working in $D$, we see that $[K_2,K_4]=1$ and $\la K_2,K_3\ra \cong Sp_4(q)$.
Now $K_0$ and $t_1$ lie in $C_G(C)$, which is an $SL_3(q)$ 
generated by long root groups in $G$.
Arguing as in Proposition \ref{e678} for $E_6(q)$, 
we deduce that $t_1$ is a root involution. If 
$K_1$ is the $SL_2(q)$ factor of its centralizer, 
then $K_1$ centralizes 
$K_3$ and $K_4$; also $\la K_1,K_2\ra = C_G(C) \cong SL_3(q)$.~\hal

\subsection{$F_4(q)$, $q$ even}\label{f4p2sec}

Assume $G$ is isomorphic to $F_4(q)$, where $q$ is even and $q>2$. 
We present an algorithm to construct basic $SL_2$ subgroups in $G$.
Recall that $\o$ denotes a generator of the multiplicative group of $\F_q$.

\begin{enumerate}
\item[1-6.]
These steps are as for the $E_6$ algorithm
in Section \ref{e6p2sec} with the following modifications. 
In Step 1, we find an element $y$ of order $(q+1)(q^3-1)$ and 
define $x = y^{q^3-1}$; 
the basic $SL_2$ subgroups constructed in Step 4 are $K_0$ and $K_1$;
in Step 6, we construct 
$D = C_G(K_0)$ which is isomorphic to $Sp_6(q)$.

\item[7.]
Construct an isomorphism $\psi$ from $D$ to the 
standard copy of $Sp_6(q)= Sp(V)$.

\item[8.]
Consider $v_2 \in K_1$. This element acts on $D$. 
Compute $T \in GL_6(q)$ such that 
$(g^{v_2})\psi = (g\psi)^T$ for all $g \in D$.

\item[9.]
Diagonalise $T$ to find a basis of $V$ with respect to which 
$T = (\o I_3,\o^{-1}I_3)$.
Choose a basis $e_1,e_2 ,e_3$ for the $\o$-eigenspace, and 
a basis $f_1,f_2,f_3$ for the $\o^{-1}$-eigenspace such 
that $(e_i,f_j) = \d_{ij}$.

\item[11.]
Define $a_0,b_0,c_0 \in D$ to be the inverse images under 
$\psi$ of the elements 
in $Sp(V)$ fixing $e_2,e_3,f_2,f_3$ and acting on $e_1,f_1$ as follows:
\[
\begin{array}{l}
a_0 \psi \,:\,e_1 \rightarrow\o^{-1} e_1,\,f_1 \rightarrow\o f_1, \\
b_0\psi \,:\, e_1 \rightarrow e_1+f_1,\,f_1 \rightarrow f_1,\\
c_0\psi \,:\, e_1 \rightarrow e_1,\,f_1 \rightarrow e_1+f_1.
\end{array}
\]
For $i=1,2$ let $a_i,b_i,c_i \in D$ be the inverse images under $\psi$ of 
the matrices in $Sp(V)$  
fixing $e_j,f_j$ for $j\ne i,i+1$ and such that 
\[
\begin{array}{l}
a_i \psi \,:\,e_i \rightarrow\o^{-1} e_i, e_{i+1} \rightarrow \o e_{i+1}, 
f_i \rightarrow\o f_i, f_{i+1} \rightarrow \o^{-1} f_{i+1}, \\
b_i \psi \,:\, e_i \rightarrow e_i+e_{i+1}, e_{i+1}\rightarrow e_{i+1}, 
f_i \rightarrow f_i, f_{i+1}\rightarrow f_{i+1}+f_i, \\
c_i \psi \,:\, e_i \rightarrow e_i, e_{i+1}\rightarrow e_i+e_{i+1}, 
f_i \rightarrow f_i+f_{i+1}, f_{i+1}\rightarrow f_{i+1}.
\end{array}
\]
Define $K_2 = \la a_0,b_0,c_0\ra$, and for $i=3,4$, 
define $K_i = \la a_{i-2},b_{i-2},c_{i-2}\ra$.
Then $K_2, K_3$ and $K_4$ are basic 
$SL_2$ subgroups in $D \cong Sp_6(q)$.
\end{enumerate}

We now have the four basic $SL_2$ subgroups 
$K_1,\ldots K_4$ in the $F_4$ Dynkin diagram.

\begin{prop}\label{f4propp2} 
The above algorithm for $F_4(q)$ for even $q$ is Monte Carlo
and runs in polynomial time.
\end{prop}
The proof is similar to that of Proposition \ref{e678propp2}.

\section{Basic $SL_2$ subgroups in $G_2(q)$}\label{g2sec}
\subsection{$G_2(q)$, $q$ odd}
Let $G$ be isomorphic to 
$G_2(q)$ with $q$ odd. 
We present an algorithm to construct basic $SL_2$ subgroups in $G$.

\begin{enumerate}
\item[1.] Find an involution $t_0 \in G$ and compute its 
centralizer $C_G(t_0) \cong (SL_2(q) \circ SL_2(q)).2$. 
Construct $S_1$ and $S_2$, the two $SL_2(q)$ factors.

\item[2.] If $q \equiv 1 \bmod 4$, then find an involution $t_1\ne t_0$ with 
$t_1 \in C_G(t_0)' = S_1S_2$; if $q \equiv 3 \bmod 4$, then  
find an involution $t_1 \in C_G(t_0) \backslash S_1S_2$. 

\item[3.] Construct the two $SL_2(q)$ factors of $C_G(t_1)$. 
For one of them -- call it $S$ -- either

\begin{enumerate}
\item[(a)] $\la S_1,S\ra \cong SL_3(q)$, $\la S_2,S \ra =G$, or

\item[(b)] $\la S_2,S\ra \cong SL_3(q)$, $\la S_1,S \ra =G$.
\end{enumerate}
\end{enumerate}

Assume (a) holds; relabel as $K_0 = S_1$, $K_1 = S$, $K_2 = S_2$. 
Now $K_1$ and $K_2$ are basic $SL_2$ subgroups, 
and we can place $K_0$, $K_1$, and $K_2$ in the 
extended $G_2$ Dynkin diagram as follows:
\[
\begin{array}{lllll}
K_0& \cdots & K_1 & \equiv > \equiv & K_2
\end{array}
\]

\begin{prop}\label{g2prop}
The above algorithm for $G_2(q)$ for odd $q$ is 
Monte Carlo and runs in polynomial time.
\end{prop}

\pf Finding the involutions and 
centralizers in Steps 1 and 2 is justified as usual using 
Lemma \ref{parwil} and Section \ref{centralizer}.
We next prove the claim in Step 3. First we show that conclusion 
(a) or (b) in that step holds for at least one involution in $C_G(t_0)$
satisfying the condition in Step 2 on being inside or outside 
the derived group.
Let $\a_1,\a_2$ be fundamental roots with $\a_1$ long, 
and let $\a_0 = 2\a_1+3\a_2$ be the highest root.
We choose notation so that, 
interchanging $S_1$ and $S_2$ if necessary, 
$S_1 = \la U_{\pm \a_0} \ra$, $S_2 = \la U_{\pm \a_2}\ra$. 
Let $t_1$ be the involution in the centre
of $\la U_{\pm \a_1} \ra$. 
If $q \equiv 1 \bmod 4$ then $t_1 \in C_G(t_0)'$ since it equals 
$h_0(i)h_2(-i) \in S_1S_2$. 
If $q \equiv 3 \bmod 4$ then $t_1 \not\in C_G(t_0)'$:
indeed, Bruhat decomposition implies that if 
$t_1 = s_1s_2$ with $s_i\in S_i$, then 
$s_1$ is in $B$ or $Bn_{\a_0}B$, and $s_2$ is in $B$ or $Bn_{\a_2}B$.  
Since $t_1 = h_1(-1) \in B$, the only possibility 
is that $s_1,s_2 \in B$, which leads to $h_1(-1) = 
u_{\a_0}(a)u_{\a_2}(a')h_0(b)h_2(b')$ for some $a,a',b,b' \in \F_q$.
This is impossible since the only involution of the form 
$h_0(b)h_2(b')$ is $t_0 = h_2(-1)$. 
When $q \equiv 1 \bmod 4$, all non-central involutions in 
$C_G(t_0)' = S_1S_2$ are conjugate; when $q \equiv 3 \bmod 4$, 
all outer involutions in $C_G(t_0)\backslash C_G(t_0)'$ are conjugate.
The claim in Step 3 follows.
\hal 

\subsection{$G_2(q)$, $q$ even}
Assume $G$ is isomorphic to $G_2(q)$, where 
$q$ is even and $q>2$. 
We present an algorithm to construct basic $SL_2$ subgroups in $G$.
Recall that $\o$ denotes a generator of the multiplicative group of $\F_q$.

\begin{enumerate}
\item[1.]
Find $y \in G$ of order $3(q-\e)$, where 
$\e = \pm 1$ and $q \not \equiv \e \bmod 3$.
Define $x = y^{q-\e}$, an element of order 3. 

\item[2.] 
If, after $O(1)$ random selections, we fail to find 
$g \in G$ with the property that $\la x,x^g \ra = K_2\cong SL_2(q)$ 
then go to Step 1.

\item[3.]
Construct an isomorphism $\phi$ from $K_2$ to the standard copy 
of $SL_2(q)$. In $K_2$, write down
\[
u = \phi^{-1} \pmatrix{1&1 \cr 0&1}, \;v = \phi^{-1} 
\pmatrix{\o^{-1} & 0 \cr 0 & \o}.
\]

\item[4.]
If $q=4$ then compute $K_0:= C_G(K_2)$. Otherwise, 
construct the involution centralizer $C_G(u)$, and 
$N: = \la C_G(u),v\ra$; apply Lemma \ref{Formula} 
to $N$ to construct $K_0 = C_N(v)'$, a long $SL_2$ subgroup 
centralizing $K_2$; now $K_0K_2 = A_1(q)\tilde A_1(q)$ in $G$.

\item[5.]
Construct an isomorphism from $K_0$ to the standard copy of $SL_2(q)$, and 
hence write down an involution $t \in K_0$. 

\item[6.]
Find $g \in G$ such that $\la K_0,t^g \ra = Y \cong SL_3(q)$.

\item[7.]
Construct an isomorphism $\phi$ from $Y$ to $SL_3(q) = SL(V)$. 
Compute $\la v_1,v_2 \ra = [V,K_0\phi]$ and $\la v_3 \ra = 
C_V(K_0\phi)$. Construct $K_1 \cong SL_2(q)$ in 
$Y$ generated by the preimages under $\phi$ of generators
for $K_1$: with respect to the basis $v_1,v_2,v_3$, 
these are
\[
\pmatrix{1&0&0 \cr 0&\o^{-1} & 0\cr 0&0&\o}, \; 
\pmatrix{1&0&0 \cr 0&1 & 1\cr 0&0&1},\;
 \pmatrix{1&0&0 \cr 0&1 & 0\cr 0&1&1}.
\]
\end{enumerate}

We now have the three $SL_2(q)$ subgroups $K_0,K_1,K_2$ 
in the extended $G_2$ Dynkin diagram.

\begin{prop}\label{g2propp2}
The above algorithm for $G_2(q)$ for even $q$ is Monte Carlo and 
runs in polynomial time.
\end{prop}

\pf In Step 1, $y$ lies in a maximal torus of $G$ contained 
in a subsystem subgroup $A_1(q)\tilde A_1(q)$, where the first factor is 
generated by long root subgroups of $G$ and the second by short root 
subgroups. Hence the element $x = y^{q-\e}$ of order 3 lies in 
$A_1(q)$ or $\tilde A_1(q)$. 
If $x \in \tilde A_1(q)$ then, 
by Proposition \ref{g2probprop}(ii),
there is a positive lower bound independent of $q$ for the probability that 
$\la x,x^g \ra$ is a conjugate of $\tilde A_1(q)$. If $x \in A_1(q)$ then 
Proposition \ref{g2probprop}(iii) shows that the probability that 
$\la x,x^g \ra$ is a conjugate of $A_1(q)$ is very small, justifying Step 2.

Consider Step 4. The construction of $C_G(u)$ is 
justified as in \cite[Theorem 3.9]{LO}, and the structure of 
this involution centralizer is 
given by \cite{AS}: $C_G(u) = QK_0$, where $Q$ is abelian of order $q^3$. 
If $U = C_{K_2}(u)$, then $C_G(u) = C_G(U)$ and 
$N_G(U) = \la C_G(U), v \ra \cong Q(K_0 \times \la v \ra ))$. 
The element $v$ acts fixed point freely on $Q$ for $q>4$, 
so we can apply Lemma \ref{Formula}. 

Finally, Proposition \ref{g2probprop} justifies Step 6: 
for random $g \in G$, there is a 
positive lower bound independent of $q$  for the probability that 
$\la K_0,t^g \ra = Y \cong SL_3(q)$. \hal

\section{Basic $SL_2$ subgroups in $^2\!E_6(q)$}\label{2e6sec}
\subsection{$^2\!E_6(q)$, $q$ odd}
Let $G$ be isomorphic to the 
quasisimple group $G(q) = \,^2\!E_6(q)$ with $q$ odd. 

Since $G(q)$ is a twisted group, we construct basic $SL_2$
subgroups and 
root elements relative to the twisted root system, 
which is of type $F_4$ (see \cite[2.4]{GLS}). Thus we aim 
to find $SL_2$ subgroups $K_1,\ldots,K_4$ forming the diagram
\[
\begin{array}{lllllll}
K_1 & - & K_2 & =>= & K_3 & - & K_4
\end{array}
\]
where $K_i \cong SL_2(q)$ for $i=1,2$; $K_i \cong SL_2(q^2)$ 
for $i=3,4$; $\la K_1,K_2 \ra \cong SL_3(q)$; 
$\la K_2,K_3\ra \cong SU_4(q)$; and 
$\la K_3,K_4\ra \cong SL_3(q^2)$ or $PSL_3(q^2)$.

We present an algorithm to construct basic $SL_2$ subgroups in $G$.
\begin{enumerate}
\item[1.] 
Find an involution $t_0 \in G$ such that $C_G(t_0)$ is of 
type $A_1\,{}^2\!A_5$. 
Construct the factors
 $K_0\cong SL_2(q)$ and $D\cong {}^2\!A_5(q)$ of the centralizer.

\item[2.]
Find an involution $t_1 \in C_G(t_0)$ such that 
$C := C_D(t_1)' \cong SL_3(q^2)$ or $PSL_3(q^2)$.

\item[3.] 
Using Section \ref{sl3q}, construct basic 
$SL_2$ subgroups $K_3,K_4$ in $C$.
 Let $t_3,t_4$ be the involutions in $K_3,K_4$. 

\item[4.]
Let $t_2 = t_0t_4$, a root involution; 
construct $K_2$, the $SL_2(q)$ factor in $C_G(t_2)$.

\item[5.]
Also $t_1$ is a root involution; construct $K_1$, 
the $SL_2(q)$ factor in $C_G(t_1)$.
\end{enumerate}

We now have the four basic $SL_2$ subgroups 
$K_1,\ldots,K_4$ in the Dynkin diagram. 

\begin{prop}\label{2e6prop}
The above algorithm for $^2\!E_6(q)$ for odd $q$ is Monte Carlo and 
runs in polynomial time.
\end{prop}

The proof is similar to that of Proposition \ref{f4prop}.

\subsection{$^2\!E_6(q)$, $q$ even}\label{2e6p2sec}
Assume $G$ is isomorphic to $^2\!E_6(q)$, where $q$ is even and $q>2$. We
present an algorithm to construct basic $SL_2$ subgroups in $G$.
Recall that $\o$ denotes a generator of the multiplicative group of $\F_q$.

\begin{enumerate}
\item[1-6.]
These steps are as for the $E_6$ algorithm in 
Section \ref{e6p2sec}, with the following modifications. 
In Step 1, we find an element $y$ of order $(q^6-1)/(3,q+1)$ and 
define $x = y^{(q^6-1)/(q+1)}$; 
the long $SL_2$ subgroups constructed 
in Step 4 are $K_0$ and $K_1$; in Step 6, we 
construct $D = C_G(K_0)$ which is isomorphic to $SU_6(q)/Z$, 
where $Z$ is a central subgroup of order 1 or $(3,q+1)$.

\item[7.] Construct an isomorphism $\psi$ from $D$ to the standard copy of
 $SU_6(q)$ modulo a central subgroup.

\item[8.] Consider $v_2 \in K_1$. This element acts
 on $D$. Compute $T \in GL_6(q^2)$ such
 that $(g^{v_2})\psi = (g\psi)^T$ for all $g \in D$. 

\item[9.] Diagonalise $T$ to find a basis of $V$ with respect to
 which $T = (\o I_3,\o^{-1}I_3)$. Choose a basis $e_1, e_2, e_3$ for
 the $\o$-eigenspace, and a basis $f_1,f_2,f_3$ for the
 $\o^{-1}$-eigenspace such that $(e_i,f_j) = \d_{ij}$.

\item[10.] Now define three basic $SL_2$ subgroups in $D$ as follows.
 Define $a_0,b_0,c_0 \in D$ to be the inverse images under 
$\psi$ of the elements in
 $SU(V)$ fixing $e_2,e_3,f_2,f_3$ and acting on $e_1,f_1$ as follows:
\[
\begin{array}{l}
a_0 \psi \,:\,e_1 \rightarrow\o^{-1} e_1,\,f_1 \rightarrow\o f_1, \\
b_0\psi \,:\, e_1 \rightarrow e_1+f_1,\,f_1 \rightarrow f_1,\\
c_0\psi \,:\, e_1 \rightarrow e_1,\,f_1 \rightarrow e_1+f_1.
\end{array}
\]
For $i=1,2$ let
$a_i,b_i,c_i \in D$ be the inverse images under $\psi$ of the 
matrices in $SU (V)$
fixing $e_j,f_j$ for $j\ne i,i+1$ and such that
\[
\begin{array}{l}
a_i \psi \,:\,e_i \rightarrow\nu^{-1} e_i, e_{i+1} \rightarrow \nu e_{i+1}, 
f_i \rightarrow\bar \nu f_i, f_{i+1} \rightarrow 
\bar \nu^{-1} f_{i+1}, \\
b_i \psi \,:\, e_i \rightarrow e_i+e_{i+1}, e_{i+1}\rightarrow e_{i+1}, 
f_i \rightarrow f_i, f_{i+1}\rightarrow f_{i+1}+f_i, \\
c_i \psi \,:\, e_i \rightarrow e_i, e_{i+1}\rightarrow e_i+e_{i+1}, 
f_i \rightarrow f_i+f_{i+1}, f_{i+1}\rightarrow f_{i+1},
\end{array}
\]
where $\nu$ is a primitive element of $\F_{q^2}$ and $\bar \nu =
\nu^q$. Define $K_2 = \la a_0,b_0,c_0\ra \cong SL_2(q)$, and for
$i=3,4$, define $K_i = \la a_{i-2},b_{i-2},c_{i-2}\ra \cong
SL_2(q^2)$. Then $K_2,K_3$ and $K_4$ are basic $SL_2$ subgroups in $D$.
\end{enumerate}

We now have the four basic $SL_2$ subgroups 
$K_1,\ldots, K_4$ in the Dynkin diagram.
\begin{prop}\label{2e6propp2} 
The above algorithm for $^2\!E_6(q)$ for even $q$ is Monte Carlo
and runs in polynomial time.
\end{prop}
The proof is similar to that of Proposition \ref{e678propp2}.

\section{Basic $SL_2$ subgroups in $^3\!D_4(q)$}\label{3d4sec}
\subsection{$^3\!D_4(q)$, $q$ odd}
Let $G$ be isomorphic to
$^3\!D_4(q)$ with $q$ odd. The twisted root system is of type $G_2$,
and we must construct basic subgroups $SL_2(q)$ and $SL_2(q^3)$.

\begin{enumerate}
\item[1.] Find an involution $t_0 \in G$ and compute its centralizer
 $C_G(t_0) \cong (SL_2(q) \circ SL_2(q^3)).2$. Construct the two
 $SL_2$ factors $K_0\cong SL_2(q)$ and $K_2\cong SL_2(q^3)$.

\item[2.] If $q \equiv 1 \hbox{ mod }4$, then find an involution 
$t_1\ne t_0$ with $t_1 \in
 C_G(t_0)' = K_0K_2$; if $q \equiv 3 \hbox{ mod }4$, then find an
 involution $t_1 \in C_G(t_0) \backslash K_0K_2$.

\item[3.] Construct the factor $K_1 \cong SL_2(q)$ of $C_G(t_1)$. 
\end{enumerate}

We now have the three basic $SL_2$ subgroups $K_0,K_1, K_2$ in 
the extended $G_2$ Dynkin diagram.

\begin{prop}\label{3d4prop}
The above algorithm for $^3\!D_4(q)$ for odd $q$ is Monte Carlo and runs
in polynomial time.
\end{prop}
The proof is similar to that of Proposition \ref{g2prop}.

\subsection{$^3\!D_4(q)$, $q$ even}\label{3d4sub}

Let $G$ be isomorphic to ${}^3\!D_4(q)$ with $q$ even.
This case differs from the others: our algorithm to construct
basic $SL_2$ subgroups employs an $O(q)$ search for an involution. 

\begin{enumerate}
\item[1.] Find an element of even order in $G$ that powers to a
root involution $t \in G$. 

\item[2.] 
Find $y \in G$ of order $(q+1)(q^3-1)$, and let $x = y^{q^3-1}$.

\item[3.]
Find $g \in G$ such that $Y:=\la x,t^g\ra \cong SL_3(q)$.

\item[4.] Construct an isomorphism from $Y$ to the standard copy of
$SL_3(q)$, and hence write down $K_0$ and $K_1$,
basic $SL_2$ subgroups in $Y$.
In $K_0$ write down the preimages of 
\[
u = \pmatrix{1&1\cr 0&1}, v = \pmatrix{\o^{-1}&0 \cr 0&\o}
\]
where $\o$ denotes a generator of the multiplicative group of $\F_q$.

\item[5.]
Construct $N := \la C_G(u), v\ra$. 

\item[6.]
For $q>2$, use Lemma \ref{Formula}
to construct $K_2:= C_N(v)'' \cong SL_2(q^3)$. 
For $q=2$ construct $K_2:=C_G(K_0)$.
\end{enumerate}

We now have the three basic $SL_2$ subgroups 
$K_0, K_1, K_2$ in the extended $G_2$ Dynkin diagram.

\begin{prop}\label{3d4propp2} 
The above algorithm for $^3\!D_4(q)$ for even $q$ is Monte Carlo
and has complexity $O(q)$.
\end{prop}

\pf By \cite[Theorem 3.8]{LO}, the proportion of 
 elements of even order in $G$ that power to a root involution 
is at least $1/8q$.
 In Step 2, $y$ lies in a subgroup $SL_2(q) \times SL_2(q^3)$, so
$x$ lies in the $SL_2(q)$ factor. In Step 3, by Proposition
\ref{3d4probprop}, there is a 
positive lower bound independent of $q$  for the probability that 
$Y:=\la x,t^g\ra \cong SL_3(q)$, a subsystem group. Step 5 yields 
$N = \la C_G(u), v\ra \cong Q.(SL_2(q^3)\times \la v \ra)$, where $Q\cong
q^{1+8}$ and $v$ acts fixed point freely on $Q/Z(Q)$. Hence 
Lemma \ref{Formula} can be applied in Step 6. 
Since $K_2$ centralizes $K_0$, this completes the extended
Dynkin diagram. \hal

\section{Labelling root and toral elements}\label{lab}
Assume that $G$ is described by a collection of generators in $GL_d(F)$,
and $G$ is isomorphic to a quasisimple exceptional group
of type $G(q)$, where $F$ and $\F_q$ are finite fields
of the same characteristic, and $q>2$. Assume also that $G(q)$ is neither
a Suzuki nor a Ree group. In previous sections we have shown how to
construct a family of basic $SL_2$ subgroups $K_r$ of $G$ as in
the Dynkin diagram. We now show how to label root and
toral elements in these subgroups consistently: we define root
elements $x_{\pm r}(c_i)$ and toral elements $h_r(\o)$ in each $K_r$,
where $c_i$ runs over an  $\F_p$-basis of $\F_q$ or an extension
field, and $\o$ is a primitive element of the field. Our
labelling algorithms are largely independent of the characteristic $p$.

We use these root and toral elements in Section \ref{highwt} to
compute the high weight of the representation of $G$ on $V = V_d(F)$ 
when $V$ is irreducible,
and in Section \ref{pres} to construct standard generators of $G$.

We summarise the result of this section.
\begin{prop}\label{labresult} 
Let $G$ be a subgroup of $GL_d(F)$,  
where $F$ is a finite field of the same
characteristic as $\F_q$, and assume  
that $G \cong G(q)$, a quasisimple group of exceptional Lie type
over $\F_q$ for $q > 2$, excluding  Suzuki and Ree groups. 
Assume also that generators are given for a family of 
basic $SL_2$ subgroups of $G$ as in
the Dynkin diagram. Subject to the availability of 
a discrete log oracle, there is a Las Vegas polynomial-time algorithm 
to label root and toral elements in each of the 
basic $SL_2$ subgroups. 
\end{prop}

The algorithm is described and justified in the remainder of this section.
We make frequent use of the algorithms to construct isomorphisms to 
various low-dimensional classical groups given by Theorem \ref{available-crec}. 

\subsection{Labelling $E_6(q)$, $E_7(q)$ and $E_8(q)$}\label{e678labsec}
Here we assume that $G\cong G(q) = E_l(q)$, $l=6, 7$ or 8. 
In Sections \ref{e678sec} and \ref{e678p2sec}, we constructed 
basic $SL_2$ subgroups $K_1,\ldots, K_l$ of $G$.

\begin{enumerate}
\item[1.] Construct an isomorphism $\phi$ from $\la K_1,K_3\ra$ to
 $SL_3(q) = SL(V)$. Choose a basis $v_1,v_2,v_3$ of $V$ such that
 $v_1 \in C_V(K_3)$, $v_2 \in [V,K_1] \cap [V,K_3]$ and $v_3 \in
 C_V(K_1)$. 
Write all matrices with respect to this basis.
 Define
\[
\begin{array}{ll}
x_1(c_i) = \phi^{-1} \pmatrix{1&c_i&0 \cr 0&1&0 \cr 0&0&1},\; 
& x_3(c_i) = \phi^{-1} \pmatrix{1&0&0 \cr 0&1&c_i \cr 0&0&1}, \\
x_{-1}(c_i) = \phi^{-1} \pmatrix{1&0&0 \cr c_i&1&0 \cr 0&0&1},\; 
& x_{-3}(c_i) = \phi^{-1} \pmatrix{1&0&0 \cr 0&1&0 \cr 0&c_i&1}, 
\end{array}
\]
where $c_i$ runs over an $\F_p$-basis of $\F_q$, and let
\[
h_1(\o) = \phi^{-1} (\o^{-1},\o,1),\; h_3(\o) = \phi^{-1} (1, \o^{-1},\o).
\]

\item[2.] Construct an isomorphism $\psi$ from $\la K_3,K_4 \ra$ to
 $SL_3(q) = SL(V)$. Choose a basis $v_1,v_2,v_3$ as in Step 1.
Compute $\l$ such that $h_3(\o) =
 (\l^{-1},\l,1)$ and define $h_4(\o) = (1,\l^{-1},\l)$. Compute
 $\mu_i$ such that
\[
x_{\pm 3}(c_i)\psi = \pmatrix{1&\mu_i&0 \cr 0&1&0 \cr 0&0&1}, 
\pmatrix{1&0&0 \cr \mu_i&1&0 \cr 0&0&1},
\]
and define
\[
x_{\pm 4}(c_i) = \psi^{-1} \pmatrix{1&0&0 \cr 0&1&\mu_i \cr 0&0&1}, 
\psi^{-1} \pmatrix{1&0&0 \cr 0&1&0 \cr 0&\mu_i&1},
\]
taking the plus terms to be both upper or both lower triangular, 
and similarly for the minus terms.

\item[3.]
Repeat Step 2 in turn for $\la K_i,K_j\ra$ with 
$(i,j) = (2,4)$, $(4,5)$, $(5,6)$, \ldots, $(l-1,l)$.
\end{enumerate} 

The justification for the above labelling is largely
self-evident.  In Step 2, observe that the
root elements $x_{\pm 3}(c_i)\psi$ are as claimed, since the
root groups generated by these elements are normalized by $h_3(\o)$.

\subsection{Labelling $F_4(q)$}

Here we assume that $G \cong G(q) = F_4(q)$. In Section \ref{f4sec} we
constructed basic $SL_2$ subgroups $K_1,\ldots, K_4$ of $G$.

\begin{enumerate}
\item[1.] Working in $\la K_1,K_2\ra \cong SL_3(q)$, label $x_{\pm
 1}(c_i), x_{\pm 2}(c_i)$ and $h_1(\o), h_2(\o)$ as in Step 1 of
 Section \ref{e678labsec}.

\item[2.] Construct an isomorphism $\psi$ from $\la K_2,K_3\ra$
 to $Sp_4(q) = Sp(V)$, and let $(\,,\,)$ be the
 associated symplectic form on $V$. Let $U = C_V(K_2\psi)$, $W = [V,
 K_2\psi ]$ and choose $e_1,f_1 \in V$ such that $e_1 \in
 C_W(x_{-2}(1)\psi)$, $f_1 \in C_W(x_2(1)\psi)$ and
 $(e_1,f_1)=1$. Let $X = \la e_1^{K_3\psi}\ra$, $Y = \la
 f_1^{K_3\psi}\ra$, so that $X$ and $Y$ are 2-spaces with $V = X \oplus
 Y$. Choose $e_2 \in U\cap X$, $f_2 \in U\cap Y$ such that
 $(e_2,f_2)=1$. Write 
 all matrices relative to the basis
 $e_1,e_2,f_2,f_1$ of $V$. 

\item[3.]
Compute $\l_i,\mu_i$ such that 
\[
x_2(c_i)\psi = \pmatrix{1&0&0&\l_i \cr 0&1&0&0 \cr 0&0&1&0 \cr 0&0&0&1}, 
x_{-2}(c_i)\psi = \pmatrix{1&0&0&0 \cr 0&1&0&0 \cr 0&0&1&0 \cr \mu_i &0&0&1}, 
\]
and define
\[
x_3(c_i) = \psi^{-1} \pmatrix{1&0&0&0 \cr \l_i&1&0&0 \cr 0&0&1&0 \cr 0&0&-\l_i&1}, 
x_{-3}(c_i) = \psi^{-1} \pmatrix{1&\mu_i&0&0 \cr 0&1&0&0 \cr 0&0&1&-\mu_i \cr 0 &0&0&1}.
\]
Compute $\l$ such that $h_2(\o)\psi = (\l^{-1},1,1,\l)$, and define
$h_3(\o) = \psi^{-1}(\l,\l^{-1},\l, \l^{-1})$.

\item[4.] Working in $\la K_3,K_4\ra \cong SL_3(q)$,
 label $x_{\pm 4}(c_i)$ and $h_4(\o)$ as in Step 2 of Section
 \ref{e678labsec}.
\end{enumerate}

Observe that the choice of basis in Step 2 ensures that the root
elements $x_{\pm 2}(c_i)\psi$, $x_{\pm 3}(c_i) \psi$ are as 
claimed in Step 3.

\subsection{Labelling $G_2(q)$ and $^3\!D_4(q)$} \label{g23d4labsec}
Here we assume that $G \cong G(q)$ is isomorphic
to either $G_2(q)$ or ${}^3\!D_4(q)$. In Sections
\ref{g2sec} and \ref{3d4sec} we constructed basic $SL_2$ subgroups
$K_0, K_1, K_2$ of $G$, with $K_0,K_1\cong SL_2(q)$ and $K_2\cong
SL_2(q)$ or $SL_2(q^3)$.

\begin{enumerate}
\item[1.] Working in $\la K_0,K_1\ra \cong SL_3(q)$, label $x_{\pm
 0}(c_i), x_{\pm 1}(c_i)$ and $h_0(\o), h_1(\o)$ as in Step 1 of
 Section \ref{e678labsec}.

\item[2.] Construct an isomorphism $\psi$ from $K_2$ to $SL(V) = SL_2(q)$
 or $SL_2(q^3)$. Compute $T \in SL(V)$ such that $(x^{h_1(\o)})\psi =
 (x\psi)^T$ for all $x \in K_2$. Choose a basis of $V$ consisting
of eigenvectors
 of $T$, and write all matrices relative to this basis.

\item[3.] Let $q=p^a$. For $q>3$ let $\Lambda = \{\o^{\pm p^j} :
 0\le j\le a-1\}$, a subset of $\F_q$ of size $2\log_pq$. 
For $\l \in \Lambda$ define $d_2(\l)
 = \psi^{-1}(\l^{-1},\l)$. Find $\l$ such that $h_1(\o)^2d_2(\l) \in
 K_0$ and define
\[
h_2(\o) = d_2(\l).
\]

\item[4.]
If $q>4$, then $\l = \o^{\e p^j}$ where $\e = \pm 1$. Define
\[
y_+ = \psi^{-1} \pmatrix{1&1\cr 0&1},\; y_- = \psi^{-1} \pmatrix{1&0\cr 1&1}.
\]
Find $\d = \pm$ such that $\la x_1(1)^{K_2},y_\d^{K_1},K_1 \ra <
G$. (For the opposite $\d$ this subgroup is $G$.) Now
$\d=\e$, except possibly for $q=5$ or 9. For $q=5$ or 9, 
if $\d\ne \e$, then replace $h_2(\o)$ by $h_2(-\o)$ or
$h_2(-\o)$ (respectively), $\e$ by $-\e$, and $j$ by 0 or 1 (respectively).

If $q=4$, find $\d$ as for $q > 4$; 
compute $i \in \{1,2\}$ such that $\l = \o^i$, 
and choose $j \in \{0,1\}$
such that $\d = (-1)^{i+j+1}$. 
If $q=3$, set $\d = \e = +$ and $j=0$.

\item[5.] Let $c_i$ ($1\le i \le a$) be an  $\F_p$-basis of $E$ 
 where $E = \F_q$ (or $\F_{q^3}$ for $^3\!D_4(q)$), and let $\nu$ be
 a primitive element of $E$ (so $\nu = \o$ if $E = \F_q$). If $\d
 =+$, then define $h_2(\nu) = \psi^{-1}(\nu^{-p^j}, \nu^{p^j})$ and
\[
x_2(c_i) = \psi^{-1} \pmatrix{1&c_i^{p^j}\cr 0&1},\; x_{-2}(c_i) = \psi^{-1} \pmatrix{1&0\cr c_i^{p^j}&1}.
\]
If $\d=-$, then define $h_2(\nu) = \psi^{-1}(\nu^{p^j}, \nu^{-p^j})$ and 
\[
x_2(c_i) = \psi^{-1} \pmatrix{1&0 \cr c_i^{p^j}&1},\; x_{-2}(c_i) = \psi^{-1} \pmatrix{1&c_i^{p^j}\cr 0&1}.
\]
\end{enumerate}

We now justify the above labelling. In Step 2, the
computation of the matrix $T$ involves solving linear equations of the
form $TA = BT$, where $A = (x^{h_1(\o)})\psi$, $B = x\psi$, for
generators $x$ of $K_2$; this system can be solved in polynomial time.
In Step 3, $\Lambda$ is introduced because 
the isomorphism $\psi$ could change eigenvalues of elements of
$K_2$ by a field automorphism or inversion. We define 
$h_2(\o)$ as in Step 3 to ensure that $h_1(\o)^2h_2(\o) = h_{2\a+3\b}(\o)
\in K_0 = \la x_{\pm 0}(c): c \in \F_q\ra $. 
Membership of 
$h_1(\o)^2d_2(\l)$ in $K_0$ can be decided readily: 
$K_0 = C_G(K_2)$, so we simply check whether the given element
commutes with the generators of $K_2$. For $q$ odd,
$h_1(-\o)^2h_2(-\o) \in K_0$; this causes a complication 
only for $q=5$ or 9: if $q$ is not one of these values, then $-\o \not
\in \Lambda$, so Step 3 defines $h_2(\o)$ uniquely. However, if
$q=5$ or 9 then $-\o = \o^{-1}$ or $\o^{-3}$ respectively, so we may need
to replace $h_2(\o)$ by $h_2(-\o)$, $\e$ by $-\e$, and $j$ by 0 or 1,
as indicated in Step 4. Similar observations apply for $q$ even: 
non-uniqueness occurs only for $q=4$ when $\o^2 = \o^{-1}$. 
In Step 4, the $y_{\pm}$ are elements of a root
group normalized by $h_1(\o)$ (this is the reason for Step 2), hence
can be taken as root elements $x_{\pm 2}(1)$. The choice of $\d$
ensures that $y_\d$ is $x_{+2}(1)$ rather than the negative, hence
justifying the parametrization of root elements in Step 5.

\subsection{Labelling $^2\!E_6(q)$} \label{2e6labsec}

Here we assume that $G \cong G(q) = {}^2\!E_6(q)$. In Section \ref{2e6sec} we
constructed basic $SL_2$ subgroups $K_1,\ldots, K_4$ of $G$, with
$K_1,K_2\cong SL_2(q)$ and $K_3,K_4\cong SL_2(q^2)$.

\begin{enumerate}
\item Construct an isomorphism $\phi$ from $\la K_2,K_3\ra$
 to $SU_4(q) = SU(V)$, let $(\,,\,)$ be the associated hermitian form
 on $V$, and write $\bar \a = \a^q$ for $\a \in \F_{q^2}$. Let $U =
 C_V(K_2\phi)$, $W = [V, K_2\phi ]$. Write $V \downarrow K_3\phi =
 X\oplus Y$ with $X,Y$ 2-spaces.
Choose $e_1 \in X\cap W$, $e_2
 \in X\cap U$, $f_1 \in Y\cap W$, $f_2 \in Y\cap U$ such that
 $(e_1,f_1) = (e_2,f_2) = \l$, where $\l$ is a fixed element of
 $\F_{q^2}$ such that $\l+\bar \l = 0$. Write all matrices with respect
 to the basis $e_1,e_2,f_2,f_1$ of $V$. 

\item
Let $c_i$ ($1\le i \le a$) be a
 $\F_p$-basis of $\F_q$, and extend it to a basis $d_i$ ($1\le i
 \le 2a$) of $\F_{q^2}$ over $\F_p$. Let $\o, \nu$ be primitive
 elements of $\F_q, \F_{q^2}$, respectively. For each $i$, define
\[
\begin{array}{ll}
x_2(c_i) = \phi^{-1} \pmatrix{1&0&0&c_i \cr 0&1&0&0 \cr 0&0&1&0 \cr 0&0&0&1}, 
& x_{-2}(c_i) = \phi^{-1} \pmatrix{1&0&0&0 \cr 0&1&0&0 \cr 0&0&1&0 \cr 
c_i &0&0&1}, \\
x_3(d_i) = \phi^{-1} \pmatrix{1&0&0&0 \cr d_i&1&0&0 \cr 0&0&1&0 \cr 
0&0&-\bar d_i&1}, 
& x_{-3}(d_i) = \phi^{-1} \pmatrix{1&d_i&0&0 \cr 0&1&0&0 \cr 
0&0&1&-\bar d_i \cr 0 &0&0&1},
\end{array}
\]
and set $h_2(\o) = \phi^{-1}(\o^{-1},1,1,\o)$, 
$h_3(\nu) = \phi^{-1}(\nu,\nu^{-1},\bar \nu, \bar \nu^{-1})$. 

\item Working in $\la K_1,K_2\ra \cong SL_3(q)$, label $x_{\pm
 1}(c_i)$ and $h_1(\o)$ as in Step 2 of Section \ref{e678labsec}.

\item Working in $\la K_3,K_4\ra \cong SL_3(q^2)$ or $PSL_3(q^2)$, label 
$x_{\pm 4}(d_i)$ and $h_4(\nu)$ as in Step 2 of Section \ref{e678labsec}.
\end{enumerate}

\section{Determining the high weight of a representation} \label{highwt}
Let $G$ be an absolutely irreducible subgroup of $GL_d(F)$ that 
is isomorphic to a quasisimple exceptional group $G(q)$ of Lie 
type over $\F_q$, where $F$ and $\F_{q}$ have the same characteristic.
Assume also that $G(q)$ is neither a Suzuki nor a Ree group. 
Write $V = V_d(F)$. In this section we describe a simple 
algorithm to compute the 
high weight of the absolutely irreducible $FG$-module $V$. 
That is, we compute the non-negative integers $n_r$ ($1 \le r \le l$) such that 
$V = V(\l)$, the irreducible module of high weight 
$\l = \sum_1^l n_r\l_r$, where 
$l$ is the rank of the corresponding simple algebraic group 
and $\l_r$ are the fundamental dominant weights. Unlike previous sections, 
the algorithm applies for all values of $q$ including $2$.

First consider the case where $G(q)$ is of untwisted type. 
The algorithm is the following.
Using the work of previous sections, 
construct the root and toral elements $x_{\pm r}(c_i)$, $h_r(\o)$ of $G$.
Construct the maximal unipotent subgroup $U$ generated by all the 
positive root elements $x_r(c_i)$. (For $F_4(2)$ additional
generators are required -- see Section \ref{f4e62}.) 
A consequence of \cite[Theorem 4.3(c)]{curtis}  is 
that $C_V(U)$ is a 1-dimensional space, 
spanned by a {\it maximal} vector $v$. Since the $h_r(\o)$
normalize $U$, they fix $C_G(U)$. Thus, 
for each $r \in \{1, \dots, l\}$, there exists 
$n_r \in \{0,\ldots, q -1 \} $ such that 
$\o^{n_r} \in F$ and
\[
v^{h_r(\o)} = \o^{n_r}v.
\]
These are the required integers $n_r$; to compute them, 
use a discrete log oracle in $\F_q$. 
The only ambiguity occurs when $v^{h_r(\o)}=v$, in which case $n_r$ can be 
0 or $q-1$. To distinguish between them, compute
the spin $\la v^{K_r}\ra$ of $v$ under $K_r$: if this is a 1-dimensional
(trivial) module for $K_r$, then $n_r=0$; otherwise $n_r=q-1$.

Now consider the twisted groups. 
For $^2\!E_6(q)$,
as in Section \ref{2e6labsec}, 
compute the root elements 
$x_{\pm r}(c_i)$, and also the toral elements 
$h_1(\o),h_2(\o)$, $h_3(\nu),h_4(\nu)$, where $\o$ and 
$\nu$ are primitive elements for
$\F_q$ and $\F_{q^2}$ respectively. 
Construct the maximal unipotent group $U$ generated by the positive
root elements $x_r(c_i)$, and compute $C_V(U) = \la v \ra$. 
Using a discrete log oracle, find $0 \le a,b,c,d,e,f \le q-1$ such that 
\[
v^{h_1(\o)} = \o^a v,\,v^{h_2(\o)} = \o^b v,\,v^{h_3(\nu)} 
= \nu^{c+dq}v,\,v^{h_4(\nu)} = \nu^{e+fq}v.
\]
The high weight of $V$ relative to the $E_6$ Dynkin diagram is
\[
\begin{array}{lllll}
e&c&b&d&f \\
 & &a& &
\end{array}
\]
Similarly, for $^3\!D_4(q)$, compute $0 \le a,b,c,d \le q-1$ such that 
\[
v^{h_1(\o)} = \o^a v,\,v^{h_2(\nu)} = \nu^{b+cq+dq^2}v
\]
where $\nu$ is now a primitive element for $\F_{q^3}$. 
The high weight of $V$ relative to the $D_4$ Dynkin diagram is $bacd$. 
In both twisted cases, we distinguish between the possibilities 
0 and $q-1$ as in the untwisted case.

Since the labelling of the root and toral elements is only determined up to an 
automorphism of $G$, the same is true of the high weight.

We have now justified the following result.
\begin{prop}\label{highwtalg} 
Subject to the availability of a discrete log oracle,
the above algorithm determines in polynomial time the 
high weight of the absolutely irreducible $FG$-module $V$, 
up to a twist by a field or graph
automorphism of $G$.
\end{prop}

\section{Constructing the standard generators}\label{pres}

Assume $G$ is described by a collection of generators in $GL_d(F)$, 
where $F$ is a finite field of the same characteristic
as $\F_q$, and $G$ is isomorphic to an exceptional group $G(q)$ of 
Lie type over $\F_q$.
Assume also that $G(q)$ is neither a Suzuki nor a Ree group. 
In previous sections we showed how to construct a family of 
basic $SL_2$ subgroups $K_r$ of $G$ as in the Dynkin diagram, 
and how to label root elements $x_{\pm r}(c_i)$ and 
toral elements $h_r(\o)$ in each $K_r$.

In this section, we use commutators among these root elements to construct 
additional root elements in rank 2 subsystems, guided by the 
Chevalley commutator relations \cite[5.2.2]{car}. 
The root elements constructed in $G$ correspond to 
the generators of the reduced Curtis-Steinberg-Tits presentation for 
$G(q)$ as in \cite[\S 4.2 and 6.1]{BGKLP}: namely,
the standard generators of $G$.  

We list these presentations on standard generators explicitly in 
Appendix~\ref{app}.  They are used to verify the correctness of 
the output of the algorithms: namely, the 
elements $x_{\pm r}(c_i)$ and $h_r(\o)$. 

We summarise the result of this section.

\begin{prop}\label{standcon} 
Let $G$ be a subgroup of $GL_d(F)$,  
where $F$ is a finite field of the same characteristic as $\F_q$ and $q > 2$, 
and assume that $G \cong G(q)$, a quasisimple group of exceptional Lie type
over $\F_q$ which is neither a Suzuki nor a Ree group. 
Assume also that generators are given for a family of 
basic $SL_2$ subgroups of $G$ as in
the Dynkin diagram. Subject to the availability of a discrete log oracle,  
there is a Las Vegas polynomial-time algorithm to construct the 
standard generators of $G$. 
\end{prop}
The proposition is justified in the following sections.

\subsection{Standard generators of 
$E_6(q)$, $E_7(q)$ and $E_8(q)$}\label{e678pres}

These are the most straight-forward cases. 
Let $l$ be the rank of $G(q)$ (so $l = 6,7$ or 8). 
From Proposition \ref{labresult} we know fundamental 
root elements $x_{\pm r}(c_i) \in G$ for 
$1\le r\le l$ and $c_i$ in an $\F_p$-basis of $\F_q$. 
For each edge $r,s$ in the Dynkin diagram with $r<s$, define additional 
root elements $x_{\pm rs}(c_i)$ by 
\[
x_{rs}(c_i) = [x_r(c_i), x_s(1)],\; x_{-rs}(c_i) = [x_{-r}(c_i), x_{-s}(-1)].
\]
The reduced Curtis-Steinberg-Tits presentation has generators 
$x_{\pm r}(c_i)$, $x_{\pm rs}(c_i)$ for all relevant $r,s,i$, 
the relations being the
Chevalley commutator relations among these elements, together with the
relations expressing that all generators have order $p$. This
presentation defines the simply connected group $E_l(q)$. 
We give an explicit version in Appendix \ref{e678app}.

If we require a presentation for the simple group $G/Z(G)$, 
 then an additional relation may be needed to kill the centre. This
only applies for $l=6$ or 7, as the simply connected group $E_8(q)$ is
simple. We know the toral elements $h_r(\o) \in K_r$. If
$Z(G) \ne 1$, then $q-1$ is divisible by 3 if $G(q) = E_6(q)$; or by
2 if $G(q) = E_7(q)$; and $Z(G) = \la z \ra$ where
\[
z = \left\{\begin{array}{l}
h_1(\l^2)h_3(\l)h_5(\l^2)h_6(\l), \hbox{ if }G(q) = E_6(q) \\
h_2(-1)h_5(-1)h_7(-1), \hbox{ if }G(q) = E_7(q)
\end{array}
\right.
\]
and $\l$ is a cube root of unity. Each $h_r(\o)$ can be
expressed in terms of the generators $x_{\pm r}(c_i)$ using the
expression
\[
h_r(\o) = n_r(\o^{-1})n_r(1)^{-1},
\]
where $n_r(c) : = x_r(c)x_{-r}(-c^{-1})x_r(c)$. 
Hence the relation $z=1$, where $z$
is as above, completes a presentation of the simple group $G/Z(G)$.

\subsection{Standard generators of $F_4(q)$}\label{f4pres}

Suppose $G(q) = F_4(q)$. From Proposition \ref{labresult} we know 
fundamental root elements
$x_{\pm r}(c_i) \in G$ for $1\le r\le 4$ and $c_i$ in an $\F_p$-basis of $\F_q$.
We define additional root elements as follows.
For $1\le r \le 4$, let 
\[
n_r = x_r(1)x_{-r}(-1)x_r(1).
\]
Now define 
\[
\begin{array}{ll}
x_{12}(c_i) = x_1(c_i)^{n_2},& x_{-12}(c_i) = x_{-1}(c_i)^{n_2},\\
x_{23}(c_i) = x_3(-c_i)^{n_2},& x_{-23}(c_i) = x_{-3}(-c_i)^{n_2},\\
x_{34}(c_i) = x_3(c_i)^{n_4},& x_{-34}(c_i) = x_{-3}(c_i)^{n_4},\\
x_{23^2}(c_i) = x_2(c_i)^{n_3},& x_{-23^2}(c_i) = x_{-2}(c_i)^{n_3}
\end{array}
\]
(where $23^2$ denotes the root $\a_2+2\a_3$ and so on). For the definition of $x_{\pm 23^2}(c_i)$ we have used the $F_4$ structure constants in \cite{GS}.

This defines all the root elements in rank 2 subsystems. 
The reduced Curtis-Steinberg-Tits presentation of $G$ 
defines the simple group $G \cong F_4(q)$, 
since this group is simply connected. We give an explicit
version in Appendix \ref{f4app}.

\subsection{Standard generators of $^2\!E_6(q)$}\label{2e6pres}

Suppose $G(q) = {}^2\!E_6(q)$.
This is very similar to the $F_4(q)$ case, except that for short roots we 
define root elements over $\F_{q^2}$ rather than $\F_q$. 
From Proposition \ref{labresult} we know 
fundamental root elements $x_{\pm r}(c_i)$ for $r=1,2$ and
$x_{\pm s}(d_i)$ for $s = 3,4$, where $c_i$ and $d_i$ run over bases for
$\F_q$ and $\F_{q^2}$ over $\F_p$, respectively. We define additional root
elements $x_{\pm 12}(c_i)$, $x_{\pm 23}(d_i)$, $x_{\pm 34}(d_i)$, $x_{\pm 23^2}(c_i)$
using exactly the same equations as in Section \ref{f4pres} for $F_4(q)$.

This defines all the root elements in rank 2 subsystems. We give an
explicit version of the reduced Curtis-Steinberg-Tits presentation of the
simply connected version of $G$ in Appendix \ref{2e6app}. This is a
variant of the presentation given in \cite[\S 6.1]{BGKLP}. 
To get a presentation for the simple group $G/Z(G)$, we
add the relation $z=1$, where $z = h_3(\l)h_4(\l^2)$; 
the $h_r$ are expressed in terms of
$x_{\pm r}(c)$ as in Section \ref{e678pres},
and $\l$ is a cube root of unity.

\subsection{Standard generators of $G_2(q)$}\label{g2pres}

Suppose $G(q) = G_2(q)$. From Proposition \ref{labresult} we know 
fundamental root elements $x_{\pm r}(c_i)$ in $G$ for $r=1,2$. 
It is convenient to change notation at
this point. Let $\a,\b$ be fundamental roots in the $G_2$ root system
with $\a$ long, $\b$ short, and write $x_{\pm 1}(c_i) = x_{\pm
 \a}(c_i)$, $x_{\pm 2}(c_i) = x_{\pm \b}(c_i)$. We define additional root
elements as follows. Let 
\[
n_\a = x_\a(1)x_{-\a}(-1)x_\a(1),\;\; n_\b = x_\b(1)x_{-\b}(-1)x_\b(1).
\]
Now define
\[
\begin{array}{ll}
x_{\a+\b}(c_i) = x_\b(-c_i)^{n_\a},& x_{-\a-\b}(c_i) = x_{-\b}(-c_i)^{n_\a},\\
x_{\a+2\b}(c_i) = x_{\a+\b}(c_i)^{n_\b},& x_{-\a-2\b}(c_i) = x_{-\a-\b}(c_i)^{n_\b},\\
x_{\a+3\b}(c_i) = x_\a(c_i)^{n_\b},& x_{-\a-3\b}(c_i) = x_{-\a}(c_i)^{n_\b},\\
x_{2\a+3\b}(c_i) = x_{\a+3\b}(-c_i)^{n_\a},& x_{-2\a-3\b}(c_i) = x_{-\a-3\b}(-c_i)^{n_\a}
\end{array}
\]

This defines all the root elements. 
The reduced Curtis-Steinberg-Tits presentation of $G$ 
defines the simple group $G \cong G_2(q)$, 
since this group is simply connected. We give an explicit
version in Appendix \ref{g2app}.

\subsection{Standard generators of $^3\!D_4(q)$}\label{3d4pres}

Suppose $G(q) = {}^3\!D_4(q)$.
This is similar to $G_2(q)$, except that for short roots we 
define root elements over $\F_{q^3}$ rather than $\F_q$. 
From Proposition \ref{labresult} we know 
root elements $x_{\pm 1}(c_i)$ for $c_i$ in an $\F_p$-basis of
$\F_q$, and $x_{\pm 2}(d_i)$ for $d_i$ in an $\F_p$-basis of $\F_{q^3}$. As
for $G_2(q)$ in Section \ref{g2pres}, relabel these as $x_{\pm \a}
(c_i)$, $x_{\pm \b}(d_i)$ respectively. 
We define additional root
elements $x_{\pm (\a+\b)}(d_i)$, $x_{\pm (\a+2\b)}(d_i)$, $x_{\pm (\a+3\b)}(c_i)$, $x_{\pm (2\a+3\b)}(c_i)$
using exactly the same equations as in Section \ref{g2pres} for $G_2(q)$.

This defines all the root elements. 
The reduced Curtis-Steinberg-Tits
presentation of $G$ 
defines the simple group $G \cong {}^3\!D_4(q)$, 
since this group is simply connected. 
We give an explicit version in Appendix \ref{3d4app}.

\section{Completion of proof of Theorem \ref{main}} \label{completion}

Theorem \ref{main} is an immediate consequence of the 
results summarised in Section \ref{prelim}, and of 
the algorithms presented and justified in 
Sections \ref{e678sec}--\ref{pres}.
Babai {\it et al.}\cite[Corollary 4.4]{BGKLP} prove that the
reduced Curtis-Steinberg-Tits presentation for 
a universal Chevalley group $G$ of rank 
at least 2 has
length $\Oh(\log^3|G|)$,  
so evaluation of the relations takes polynomial time.
That the resulting constructive recognition
algorithm is Las Vegas is established by
verifying that the standard generators satisfy
these presentations, which are given explicitly in Appendix \ref{app}.

\section{Algorithms for $q=2$}\label{qeq2}


Our algorithms to construct basic $SL_2$ subgroups 
fail when $q=2$: the critical elements
$v_1$ and $v_2$ constructed in Step 4 of Section \ref{e6p2sec} 
are now both the identity; 
the algorithms to construct standard generators 
in Section \ref{pres} also fail in some cases. 

Since it is desirable to have practical recognition algorithms for exceptional 
groups over $\F_2$, we now provide such.  We often exploit the fact 
that explicit computations can be performed readily in some of their 
subgroups using standard machinery; for details of such, see, for 
example, \cite[Chapter 4]{HoltEickOBrien05}. 
We omit $G_2(2)$ since it is isomorphic 
to the almost simple classical group $U_3(3).2$.

\subsection{$E_6(2), F_4(2)$ and $\,^2\!E_6(2)$}\label{f4e62}
Assume $G$ is isomorphic to one of $E_6(2), F_4(2)$ or $\,^2\!E_6(2)$.

\begin{enumerate}
\item Apply Steps 1-4 of Sections \ref{e6p2sec}, \ref{f4p2sec} and  
\ref{2e6p2sec}. 
These construct basic $SL_2$ subgroups $K_0,K_i$ where 
$i=2$ for $E_6(2)$ and $i=1$ for $F_4(2)$ and $^2\!E_6(2)$. 
They also find a root involution $u:= u_1^+ \in K_0$.

\item Construct $C_G(u) = QD$, where $D = C_G(K_0)$ 
and $Q$ is a normal 2-subgroup. 

\item Construct $Q$, the soluble radical of $C_G(u)$. 
Now construct $D$ as follows. 
Find involutions $s \in D\backslash Q$
such that $|Q:C_Q(s)| \le 2^{12}$. Note that $C_Q(s)$ can be computed 
from the action of $s$ on the vector space $Q/Z(Q)$ over $\F_2$. 
Search for sufficient $Q$-conjugates of $s$ lying 
in $C_G(K_0)$ to generate $D$.

\item Find an involution $t \in D$ such that $\la K_i,t \ra \cong SL_3(2)$. 
Compute $D \cap \la K_i,t\ra \cong SL_2(2)$, and call 
it $K_j$, where $j=4$ for $G\cong E_6(2)$ and $j=2$ otherwise.

\item Construct $T := C_D(K_i)$ which
is isomorphic to $SL_3(2)^2$, $SL_3(2)$ or $SL_3(4)$ 
for $G \cong E_6(2)$, $F_4(2)$ or $^2\!E_6(2)$ respectively. 

\item 
Compute $C_T(K_j)$  which is isomorphic to
$SL_2(2)^2$, $SL_2(2)$ or $SL_2(4)$ respectively. 
Define its $SL_2$ factors to be 
$K_1K_6$, $K_4$ or $K_4$, respectively.

\item Search in $T$ for the remaining basic $SL_2$ subgroups 
$K_3K_5 \cong SL_2(2)^2$, $K_3 \cong SL_2(2)$ or $K_3 \cong SL_2(4)$. 
We now know all the basic $SL_2$ subgroups in $G$.

\item The labelling of fundamental root elements is carried out 
as in Section \ref{lab}. 

\item The construction of the standard generators is as in Sections \ref{e678pres}, \ref{f4pres} and \ref{2e6pres}.
\end{enumerate}

\subsection{$E_7(2)$}\label{e72}
Assume $G \cong E_7(2)$.
\begin{enumerate}
\item Construct $K_0,K_1$ and $D = C_G(K_0)$ as in the previous section. 
Find an involution $t\in D$ such that 
$\la K_1,t\ra \cong SL_3(2)$, and compute $K_3 = D \cap \la K_1,t\ra$.

\item Step 5 of the previous section is too expensive to apply in 
$D \cong \O_{12}^+(2)$. Instead, we first compute 
$C_D(K_3) \cong SL_2(2) \times \O_8^+(2)$. 
Name the direct factors as $K_2$ and $E$ respectively.

\item 
Compute $C_E(K_1) \cong SL_4(2)$, and construct $K_5,K_6,K_7$, 
basic $SL_2$ subgroups in $C_E(K_1)$. 

\item In $C_D(K_6,K_7) \cong SL_4(2)$, search for the remaining 
basic $SL_2$ subgroup $K_4\cong SL_2(2)$, satisfying 
$[K_1,K_4]=1$ and $\la K_4,K_i \ra \cong SL_3(2)$ for $i=2,3,5$. 

\end{enumerate}
We now know all of the basic $SL_2$ subgroups in $G$. The labelling of 
fundamental root elements and the construction of standard generators 
is unchanged from Sections \ref{e678labsec} and \ref{e678pres}.

\subsection{$E_8(2)$}
An approach modelled on the previous algorithms is too expensive to 
apply to $E_8(2)$.  Instead, we present a different algorithm which 
recognises the 
group only in its 248-dimensional adjoint representation.

Assume that $G \leq GL_{248}(F)$, 
where $F$ is a finite field of characteristic 2, and $G \cong E_8(2)$. 
\begin{enumerate}
\item Using \cite{GLO}, find a basis of $V = V_{248}(F)$ with respect to which 
the action of $G$ is realised over $\F_2$. 
Replace $V$ by the $\F_2$-span of this basis, and $F$ by $\F_2$.

\item Construct $K_0$ and $D = C_G(K_0) \cong E_7(2)$ as in 
Steps 1-3 of Section \ref{f4e62}.

\item Apply the algorithm of Section \ref{e72} to construct basic 
$SL_2$ subgroups and 
root elements $x_{\pm r}(1)$ ($1\le r\le 7$) in $D$. 
Let $x_{\pm 0}(1)$ be two involutions generating $K_0$.

\item Let $\hat G$ be the standard copy of $E_8(2)$ in $GL_{248}(2)$, with 
fundamental root element generators $\hat x_{\pm r}(1)$ ($1\le r \le 8$). 
Let $\hat x_{\pm 0}(1)$ be root elements in $\hat G$ corresponding to the 
longest root $\a_0$ in the root system.

\item Compute all matrices $g \in GL_{248}(2)$ such that 
$x_{\pm r}(1)^g = \hat x_{\pm r}(1)$ for $0 \le r \le 7$. There are 
precisely two such matrices. To see this, observe that 
$\la x_{\pm r}(1) : 0 \le r \le 7 \ra = DK_0$, so any two such matrices 
$g$ differ by an element of $C_{GL(V)}(DK_0)$. 
Now $V\downarrow DK_0 = W_1\oplus W_2\oplus W_3$, a sum of indecomposables of 
dimensions 2, 112 and 134. Here $W_1$ and $W_2$ are irreducible and $W_3$ is 
uniserial with socle series having irreducible factors of 
dimensions $1, 132, 1$. 
Let $M$ be the maximal submodule of $W_3$, and 
let $\la s \ra$ be the socle of $M$. 
Then ${\rm Hom}_{DK_0}(V,V)$ has dimension 4, with basis $1, \pi_1,\pi_2,\phi$ 
where $\pi_i$ is the projection $V \rightarrow W_i$ and $\phi$ sends 
$w_1+w_2+w_3$ to 0 if $w_3 \in M$ and to $s$ if $w_3 \not \in M$. 
The only invertible elements of this space are $1$ and $1+\phi$. 
Hence, as asserted, there are exactly two such matrices $g$. 
Call them $g_1,g_2$.

\item For $i=1,2$ define $x_{\pm 8}^{(i)} = \hat x_{\pm 8}(1)^{g_i^{-1}}$. 
Decide for which value of $i$ the 
group $\la D, x_{\pm 8}^{(i)} \ra$ is isomorphic to $E_8(2)$;
the other is a ``large" subgroup of 
$GL_{248}(2)$ containing elements of order much larger than those of $E_8(2)$.
For this value of $i$, define $x_{\pm 8}(1) = x_{\pm 8}^{(i)}$. 
\end{enumerate}

We have now labelled all fundamental root elements 
$x_{\pm r}(1)$ ($1\le r\le 8$) in $G$. The construction of 
the standard generators is unchanged from Section \ref{pres}.

\subsection{$^3\!D_4(2)$}
Assume $G \cong \,^3\!D_4(2)$. 
Let $\a,\b$ be fundamental long and short roots in the root system, as 
in Section \ref{g2pres}.
\begin{enumerate}

\item Construct subgroups $K_0,K_1 \cong SL_2(2)$ and 
$K_2 \cong SL_2(8)$ as in Section \ref{3d4sub}.

\item Construct an isomorphism from $K_2$ to $SL_2(8)$. 
In $K_2$, define $x_\b(c), x_{-\b}(c)$ (for $c \in \F_8$) to be the 
preimages of $\pmatrix{1&c \cr 0&1}$, $\pmatrix{1&0 \cr c&1}$ respectively.

\item Let $X_{\pm \b} = \{x_{\pm \b}(c) : c \in \F_8 \}$. 
Find $g \in K_2$ and involutions $x_{\e} \in K_1$  ($\e = \pm$) such that 
each of $\la x_{\e}^{K_2}, X_{\e \b}^g, K_1\ra$ is a proper subgroup of $G$. 
Define $x_{\pm \a}(1) = x_{\pm}$ and replace 
$x_{\pm \beta} (c)$ by $x_{\pm \beta} (c)^g$.
Now we have labelled the fundamental root elements of $G$. 


\item Construct the remaining standard generators of $G$ as in Section \ref{3d4pres}.

\end{enumerate}

\section{Implementation and performance}\label{imp}

We have implemented these algorithms in {\sc Magma}. 
We use 
the product replacement algorithm \cite{Celleretal95} to generate 
random elements; our implementations of \cite{Bray}, 
\cite{Brooksbank03a}, 
\cite{CLGO}, 
\cite{DLGO14}, 
and \cite{sl3q}; and Brooksbank's implementations of 
his algorithm \cite{Brooksbank08} for 
constructive recognition of $Sp_4(q)$.

The computations reported in Table \ref{table1} were carried out
using {\sc Magma} V2.19 on a 2.8 GHz processor.
We list the CPU time $t_1$ in seconds taken to construct
standard generators in a random conjugate of the
standard copy of dimension $d_1$ of an exceptional group of type $G(q)$; 
sometimes, we list $t_2$, the time taken to perform 
the same task in an irreducible representation of dimension $d_2$. 
The time is averaged over three runs.

We use Taylor's implementation of \cite{CMT, CT} to write an element of $G(q)$ as a word in the standard generators. As one illustration, it takes 17 seconds to write an element of $E_8(5^2)$ as a word 
in its standard generators.

\begin{table}[ht]\label{tabRTStdGens}
\centering 
\begin{tabular}{lrrrr}
Group     & $d_1$ & $t_1$ & $d_2$ & $t_2$ \\\hline\hline

$E_6(2^3)$ & 27 & 8 & 78 & 51\\ \hline 
$E_6(5^2)$ & 27 & 15  & 78 & 119 \\ \hline 

$E_7(2^{3})$ & 56 & 35 & 133 & 158 \\\hline
$E_7(5^2)$ &56 & 53 &133 & 301 \\ \hline 

$E_8(2^3)$ & 248 & 978 & -- & -- \\     \hline 
$E_8(5^{2})$ & 248 & 520 & -- & --\\\hline

$F_4(2)$ & 26 & 11 & 246 & 235 \\ \hline 
$F_4(2^3)$ & 26 & 14 & 246 & 607 \\ \hline 
$F_4(5^3)$ & 26 & 30 & 52 & 248 \\ \hline 

$G_2(2^3)$ & 6 & 1 & 14 & 2\\ \hline 
$G_2(5^3)$ & 7 & 2 & 14 & 3\\ \hline 

${}^2E_6(2^{3})$ & 27 & 99 & 78 & 790\\\hline
${}^2E_6(5^{2})$ & 27 & 102 & 78 & 865\\\hline

${}^3D_4(2^{6})$ & 8 & 15 & 26 & 122\\\hline
${}^3D_4(5^{3})$ & 8 & 6 & 28 & 90\\\hline

\end{tabular}
\vspace{1ex}
\caption{Time to construct standard generators}\label{table1}
\end{table}

\appendix 
\section{Reduced Curtis-Steinberg-Tits presentations}\label{app}

The Curtis-Steinberg-Tits presentations are well known; the
reduced versions using only an $\F_p$-basis of the field $\F_q$ (and
extensions) are described in 
\cite{BGKLP}. Since we know of no explicit versions listing
the constants in the Chevalley relations, which 
we need for our work, we include such here. 
The constants are calculated using \cite[5.2.2]{car} together
with the $N_{\a\b}$ structure constants for the $G_2$ and $F_4$ Lie
algebras from \cite{GS}.

In all cases, the generators are the root elements we have constructed
in Section \ref{pres}, namely the elements $x_r(c_i)$ for roots $r$ in
subsystems spanned by two non-orthogonal fundamental roots, and
elements $c_i$ in an $\F_p$-basis of $\F_q$ (or an extension field). 
In every case, the presentation contains the following relations: 
\[
\begin{array}{ll}
 & x_r(c_i)^p = 1, \\
& [x_r(c_i),x_s(d_i)] = 1 \hbox{ if } r+s \hbox{ is not a root}.
\end{array}
\]
For $c = \sum k_ic_i \in \F_q$ (or an extension field) with $k_i \in
\F_p$, we set $x_r(c) = \prod x_r(c_i)^{k_i}$.

We present the remaining relations for each type below.

\subsection{Relations for $E_6(q),E_7(q),E_8(q)$}\label{e678app}

The relations for these types are simple: for each edge $rs$ in the
Dynkin diagram with $r<s$, and for $c_i, d_i$ in 
the $\F_p$-basis of $\F_q$,
\[
\begin{array}{ll}
(1) & [x_r(c_i),x_s(d_i)] = x_{rs}(c_id_i) \\
(2) & [x_{-r}(c_i),x_{-s}(d_i)] = x_{-rs}(-c_id_i) \\
(3) & [x_r(c_i),x_{-rs}(d_i)] = x_{-s}(-c_id_i) \\
(4) & [x_s(c_i),x_{-rs}(d_i)] = x_{-r}(c_id_i) \\
(5) & [x_{-r}(c_i),x_{rs}(d_i)] = x_{s}(c_id_i) \\
(6) & [x_{-s}(c_i),x_{rs}(d_i)] = x_{r}(-c_id_i) 
\end{array}
\]

\subsection{Relations for $F_4(q)$}\label{f4app}

For the edges 12 and 34 in the Dynkin diagram of $F_4$ we have the
relations (1)-(4) of the previous section. The remaining relations are
the following: 
\[
\begin{array}{ll}
(1) & [x_2(c_i),x_3(d_i)] = x_{23}(c_id_i)x_{23^2}(c_id_i^2) \\
(2) & [x_2(c_i),x_{-23}(d_i)] = x_{-3}(-c_id_i)x_{-23^2}(-c_id_i^2) \\
(3) & [x_{-2}(c_i),x_{23}(d_i)] = x_{3}(c_id_i)x_{23^2}(-c_id_i^2) \\
(4) & [x_{-2}(c_i),x_{-3}(d_i)] = x_{-23}(-c_id_i)x_{-23^2}(c_id_i^2) \\
 (5) & [x_{23^2}(c_i),x_{-3}(d_i)] = x_{23}(c_id_i)x_{2}(c_id_i^2) \\
(6) & [x_{23^2}(c_i),x_{-23}(d_i)] = x_{3}(-c_id_i)x_{-2}(-c_id_i^2) \\
(7) & [x_{-23^2}(c_i),x_{23}(d_i)] = x_{-3}(c_id_i)x_{2}(-c_id_i^2) \\
(8) & [x_{-23^2}(c_i),x_{3}(d_i)] = x_{-23}(-c_id_i)x_{-2}(c_id_i^2) \\
(9) & [x_{23}(c_i),x_3(d_i)] = x_{23^2}(2c_id_i) \\
(10) & [x_{-23}(c_i),x_{-3}(d_i)] = x_{-23^2}(-2c_id_i) \\
(11) & [x_{23}(c_i),x_{-3}(d_i)] = x_{2}(2c_id_i) \\
(12) & [x_{3}(c_i),x_{-23}(d_i)] = x_{-2}(2c_id_i)
\end{array}
\]

\subsection{Relations for $^2\!E_6(q)$}\label{2e6app}

Here the Dynkin diagram is $F_4$. For the edge 12 we have relations
(1)-(4) of Section \ref{e678app} with $c_i,d_i$ in an $\F_p$-basis of $\F_q$,
and for edge 34 we have these relations for $c_i,d_i$ in an $\F_p$-basis of
$\F_{q^2}$. The remaining relations are (1)-(12) in Appendix
\ref{f4app} with the following adjustments:

\begin{enumerate}
\item[(a)] in relations (1)-(8), $c_i \in \F_q$, $d_i\in
 \F_{q^2}$, and in (9)-(12) $c_i,d_i\in \F_{q^2}$;
\item[(b)] in relations (1)-(4), $c_id_i^2$ is replaced by 
$c_id_i\overline d_i$ (where $\overline d_i = d_i^q$);
\item[(c)] in relations (5)-(8), $c_id_i$ is replaced by 
$c_i\overline d_i$, and $c_id_i^2$ is replaced by $c_id_i\overline d_i$;
\item[(d)] in relations (9)-(10), $2c_id_i$ is replaced by 
$c_i\overline d_i+\overline c_id_i$;
\item[(e)] in relations (11)-(12), $2c_id_i$ is replaced by 
$c_id_i+\overline c_i\overline d_i$.
\end{enumerate}
Note that this is a variant of the presentation described in \cite[6.1]{BGKLP}. 


\subsection{Relations for $G_2(q)$}\label{g2app}

As in Section \ref{g2pres}, we let $\a,\b$ be fundamental roots in the
$G_2$ root system, and define root elements $x_r(c_i)$ for $r$ one of
the long roots $\pm \a$, $\pm (\a+3\b)$, $\pm (2\a+3\b)$, or one of
the short roots $\pm \b$, $\pm (\a+\b)$, $\pm (\a+2\b)$ and $c_i$ in a
$\F_p$-basis of $\F_q$. The relations are the following:

\[
\begin{array}{ll}
(1) & [x_\a (c_i),x_\b (d_i)] = x_{\a+\b}(c_id_i)x_{\a+2\b}(-c_id_i^2)x_{\a+3\b}(c_id_i^3)x_{2\a+3\b}(c_i^2d_i^3) \\
(2) & [x_\a (c_i),x_{-\a-\b} (d_i)] = x_{-\b}(-c_id_i)x_{-\a-2\b}(c_id_i^2)x_{-2\a-3\b}(c_id_i^3)x_{-\a-3\b}(-c_i^2d_i^3) \\
(3) & [x_{\a+3\b} (c_i),x_{-\b} (d_i)] = x_{\a+2\b}(-c_id_i)x_{\a+\b}(-c_id_i^2)x_{\a}(c_id_i^3)x_{2\a+3\b}(-c_i^2d_i^3) \\
(4) & [x_{\a+3\b} (c_i),x_{-\a-2\b} (d_i)] = x_{\b}(c_id_i)x_{-\a-\b}(-c_id_i^2)x_{-2\a-3\b}(c_id_i^3)x_{-\a}(c_i^2d_i^3) \\
(5) & [x_{2\a+3\b} (c_i),x_{-\a-\b} (d_i)] = x_{\a+2\b}(-c_id_i)x_{\b}(-c_id_i^2)x_{-\a}(-c_id_i^3)x_{\a+3\b}(c_i^2d_i^3) \\
(6) & [x_{2\a+3\b} (c_i),x_{-\a-2\b} (d_i)] = x_{\a+\b}(c_id_i)x_{-\b}(c_id_i^2)x_{-\a-3\b}(-c_id_i^3)x_{\a}(-c_i^2d_i^3) \\
(7) & [x_{-\a} (c_i),x_{-\b} (d_i)] = x_{-\a-\b}(-c_id_i)x_{-\a-2\b}(-c_id_i^2)x_{-\a-3\b}(-c_id_i^3)x_{-2\a-3\b}(c_i^2d_i^3) \\
(8) & [x_{-\a} (c_i),x_{\a+\b} (d_i)] = x_{\b}(c_id_i)x_{\a+2\b}(c_id_i^2)x_{2\a+3\b}(-c_id_i^3)x_{\a+3\b}(-c_i^2d_i^3) \\
(9) & [x_{-\a-3\b} (c_i),x_{\b} (d_i)] = x_{-\a-2\b}(c_id_i)x_{-\a-\b}(c_id_i^2)x_{-\a}(-c_id_i^3)x_{-2\a-3\b}(-c_i^2d_i^3) \\
(10) & [x_{-\a-3\b} (c_i),x_{\a+2\b} (d_i)] = x_{-\b}(-c_id_i)x_{\a+\b}(-c_id_i^2)x_{2\a+3\b}(-c_id_i^3)x_{\a}(c_i^2d_i^3) \\
(11) & [x_{-2\a-3\b} (c_i),x_{\a+\b} (d_i)] = x_{-\a-2\b}(c_id_i)x_{-\b}(-c_id_i^2)x_{\a}(c_id_i^3)x_{-\a-3\b}(c_i^2d_i^3) \\
(12) & [x_{-2\a-3\b} (c_i),x_{\a+2\b} (d_i)] = x_{-\a-\b}(-c_id_i)x_{\b}(c_id_i^2)x_{\a+3\b}(c_id_i^3)x_{-\a}(-c_i^2d_i^3) \\
(13) & [x_{\b} (c_i),x_{\a+\b} (d_i)] = x_{\a+2\b}(2c_id_i)x_{\a+3\b}(-3c_i^2d_i)x_{2\a+3\b}(-3c_id_i^2) \\
(14) & [x_{\b} (c_i),x_{-\a-2\b} (d_i)] = x_{-\a-\b}(-2c_id_i)x_{-\a}(3c_i^2d_i)x_{-2\a-3\b}(3c_id_i^2) \\
(15) & [x_{\a+\b} (c_i),x_{-\a-2\b} (d_i)] = x_{-\b}(2c_id_i)x_{\a}(-3c_i^2d_i)x_{-\a-3\b}(-3c_id_i^2) \\
(16) & [x_{\a+2\b} (c_i),x_{-\b} (d_i)] = x_{\a+\b}(-2c_id_i)x_{2\a+2\b}(-3c_i^2d_i)x_{\a}(-3c_id_i^2) \\
(17) & [x_{\a+2\b} (c_i),x_{-\a-\b} (d_i)] = x_{\b}(2c_id_i)x_{\a+3\b}(3c_i^2d_i)x_{-\a}(3c_id_i^2) \\
(18) & [x_{-\b} (c_i),x_{-\a-\b} (d_i)] = x_{-\a-2\b}(-2c_id_i)x_{-\a-3\b}(-3c_i^2d_i)x_{-2\a-3\b}(-3c_id_i^2) \\
(19) & [x_{\a} (c_i),x_{\a+3\b} (d_i)] = x_{2\a+3\b}(c_id_i) \\
(20) & [x_{\a} (c_i),x_{-2\a-3\b} (d_i)] = x_{-\a-3\b}(-c_id_i) \\
(21) & [x_{\a+3\b} (c_i),x_{-2\a-3\b} (d_i)] = x_{-\a}(c_id_i) \\
(22) & [x_{2\a+3\b} (c_i),x_{-\a} (d_i)] = x_{\a+3\b}(-c_id_i) \\
(23) & [x_{2\a+3\b} (c_i),x_{-\a-3\b} (d_i)] = x_{\a}(c_id_i) \\
(24) & [x_{-\a} (c_i),x_{-\a-3\b} (d_i)] = x_{-2\a-3\b}(-c_id_i) \\
(25) & [x_{\b} (c_i),x_{\a+2\b} (d_i)] = x_{\a+3\b}(3c_id_i) \\
(26) & [x_{\b} (c_i),x_{-\a-\b} (d_i)] = x_{-\a}(3c_id_i) \\
(27) & [x_{\a+\b} (c_i),x_{\a+2\b} (d_i)] = x_{2\a+3\b}(3c_id_i) \\
(28) & [x_{\a+\b} (c_i),x_{-\b} (d_i)] = x_{\a}(3c_id_i) \\
(29) & [x_{-\b} (c_i),x_{-\a-2\b} (d_i)] = x_{-\a-3\b}(-3c_id_i) \\
(30) & [x_{-\a-\b} (c_i),x_{-\a-2\b} (d_i)] = x_{-2\a-3\b}(-3c_id_i) 
\end{array}
\]

 \subsection{Relations for $^3\!D_4(q)$}\label{3d4app}
Here the Dynkin diagram is $G_2$. The relations are (1)-(30) in
Appendix \ref{g2app} with the following adjustments:

\begin{enumerate}
\item[(a)] in relations (1)-(12), 
$c_i \in \F_q$, $d_i\in \F_{q^3}$; 
in (13)-(18), $c_i,d_i\in \F_{q^3}$; 
in (19)-(24), $c_i,d_i \in \F_q$; and 
in (25)-(30), $c_i,d_i \in \F_{q^3}$; 
\item[(b)] in relations (1)-(12), $c_id_i^2$ is 
replaced by $c_i\bar d_i\bar{\bar {d_i}}$ 
(where $\bar d_i = d_i^q$), $c_id_i^3$ by 
$c_id_i\bar d_i \bar{\bar {d_i}}$, and 
$c_i^2d_i^3$ by $c_i^2d_i\bar d_i \bar{\bar {d_i}}$; 
\item[(c)] in relations (13)-(18), $2c_id_i$ is replaced by 
$\bar c_i\bar{\bar {d_i}}+ \bar{\bar{c_i}}\bar d_i$, 
$3c_i^2d_i$ by 
$c_i\bar c_i\bar{\bar{d_i}}+\bar c_i \bar{\bar {c_i}}d_i+\bar{\bar{c_i}} 
c_i\bar d_i$, 
and $3c_id_i^2$ by $c_i\bar d_i\bar{\bar{d_i}}+\bar 
c_i \bar{\bar {d_i}}d_i+\bar{\bar{c_i}} d_i\bar d_i$; 
\item[(d)] in relations (25)-(30), $3c_id_i$ is replaced 
by $c_i d_i+\bar c_i\bar d_i+ \bar{\bar{c_i}}\bar{\bar{d_i}}$.
\end{enumerate}

\end{document}